# Generalized pair-wise logit dynamic and its connection to a mean field game: theoretical and computational investigations focusing on resource management

Hidekazu Yoshioka[1, *], Motoh Tsujimura[2]

[1] Japan Advanced Institute of Science and Technology, 1-1 Asahidai, Nomi, Ishikawa 923-1292, Japan

[2] Doshisha University, Karasuma-Higashi-iru, Imadegawa-dori, Kamigyo-ku, Kyoto 602-8580, Japan

[*] Corresponding author: yoshih@jaist.ac.jp

ORCID:0000-0002-5293-3246 (Hidekazu Yoshioka), 0000-0001-6975-9304 (Motoh Tsujimura)

**Abstract** Logit dynamics are evolution equations that describe transitions to equilibria of actions among many players. We formulate a pair-wise logit dynamic in a continuous action space with a generalized exponential function, which we call a generalized pair-wise logit dynamic, depicted by a new evolution equation nonlocal in space. We prove the well-posedness and approximability of the generalized pair-wise logit dynamic to show that it is computationally implementable. We also show that this dynamic has an explicit connection to a mean field game of a controlled pure-jump process, with which the two different mathematical models can be understood in a unified way. Particularly, we heuristically discuss that the generalized pair-wise logit dynamic is derived as a myopic version of the corresponding mean field game, and that the conditions to guarantee the existence of unique solutions are different from each other. The key in this procedure is to find the objective function to be optimized in the mean field game based on the logit function. The monotonicity of the utility is unnecessary for the generalized pair-wise logit dynamic but crucial for the mean field game. Finally, we present applications of the two approaches to fisheries management problems with collected data.

# 1. Introduction

## 1.1 Background

The management of environment and resources has been a global challenge in the path toward sustainable development (Chen et al., 2023; Fu et al., 2023; Ulucak and Baloch, 2023), involving decision-making among many socially interacting players, such as stakeholders, where each player has its own utility. From a dynamical system viewpoint, a common ground among the players appears as an equilibrium of a system of equations and/or inequalities. Such equilibria can be found effectively using game-theoretic approaches. For example, Rettieva (2023) studied multicriteria games and investigated the cooperate incentive equilibrium. Wang et al. (2023) focused on a static game to investigate sustainable marine resources. Ahmadi Forushani et al. (2023) formulated a signaling game for evaluating water-allocation competitiveness subject to asymmetric information.

Decision-making problems have been studied in evolutionary game theory to track transient to equilibrium behaviors of interacting players parameterized through their actions and utilities (e.g., Sandholm, 2020). A major research subject in this theory is evolution equations, such as systems of ordinary differential equations (ODEs) and/or partial differential equations (PDEs), which govern the utility-driven actions of players. Major mathematical models among them are the replicator dynamic that covers both transient and stationary equilibria (Mendoza-Palacios and Hernández-Lerma, 2015); logit dynamic in a continuous action space as an approximation to Nash equilibria under noisy environments (Lahkar and Riedel, 2015); gradient dynamic that reduces the problem objective to solving a nonlinear hyperbolic equation (Friedman and Ostrov, 2013; Lai et al., 2024); dynamics generalizing the aforementioned ones (Cheung, 2016; Zusai et al., 2023); and stochastic dynamics (Choudhury and Aydinyan, 2023). Engineering applications of evolutionary game theory include but are not limited to groundwater management (Biancardi et al., 2023), environmental taxation (Sartzetakis et al., 2023), transboundary pollution control (Huang, 2023), and interactions between social media and ecological conservation (Takashina et al., 2023). Related models such as the Ising (Drechsler, 2023a; Drechsler, 2023b) and Darwinian models also exist (Benjamin and Arjun Krishnan, 2023). The optimal action obtained from an evolutionary game model typically depends only on the current state and is therefore not anticipative.

Another stream for modeling decision-making problems among players is the mean field game (MFG) theory, where infinitely many interacting players choose their actions considering a current to future expectation (e.g., Lasry and Lions, 2007); utilities gained in the future are considered by the players through a conditional expectation (i.e., future prediction). The governing equation of the optimal actions of players in this theory is a coupled system of time-forward Fokker–Planck (FP) equation of the population measure and a time-backward Hamilton–Jacobi–Bellman (HJB) equation as its formal nonlinear adjoint. The evolution equation to be solved is more complicated in an MFG than in an evolutionary game because the former deals with a forward–backward system. A vast number of studies have been conducted to investigate the well-posedness of MFG models, focusing on the monotonicity of the utility/cost (Gomes and Saúde, 2021; Cohen and Zell, 2023), analytical solutions (Dai Pra et al., 2023; Ullmo et al., 2019), and numerical approximations (Ghilli et al., 2023; Yu et al., 2023). The MFG also finds its applications in manufacturing (Wang and Xin, 2023), epidemic control (Roy et al., 2023), crowd evacuation (Yano and Kuroda, 2023), and management of electric vehicles (Lin et al., 2022).

Recently, it has been reported that an evolutionary game model appears to be a myopic version of an MFG for a diffusion process (Barker, 2019). Subsequently, the two approaches have been systematically compared against different utility profiles (Barker et al., 2022). A key in this context is to scale the objective function in the MFG to relate it to the corresponding evolutionary game. It is reasonable to expect that a similar linkage is found under a different setting such as a jump-driven dynamic; however, to the best of the authors' knowledge, investigations in this direction have been insufficient. Such studies would contribute not only to theory of dynamic games, but also to several application areas where social behavior of players is crucial, which include but are not limited to epidemiology (Certório et al., 2022), ecology (Thon et al., 2024), and resource management (Gao et al., 2024).

### 1.2 Objective and contribution

The objectives of this study are to formulate and study a new evolutionary game model, to find its MFG counterpart, and to compare them theoretically and numerically through applications. Our contributions to achieve these objectives are explained as follows.

The evolution equation that we focus on is a pair-wise logit model (e.g., Aiba, 2015; Sandholm and Staudigl, 2018; Leon et al., 2023) extended to a continuous action space (e.g., Lahkar and Riedel, 2015). Particularly, we define a generalized logit function, i.e., the softmax function for the decision-making in noisy environments, using a $\kappa$-exponential function (e.g., Kaniadakis et al., 2020) instead of an exponential one. The $\kappa$-exponential function reduces to an exponential one in a limit case, and hence generalizes the classical pair-wise logit model. The generalization here implies considering a wider range of uncertainties in the decision-making of players. Note that seeking for generic logit functions, more specifically sigmoidal activation functions, is also important in machine learning (Apicella et al., 2021; Dombi and Jónás, 2022). Moreover, the evolutionary game theory for machine learning has recently been a popular research topic, to which this study may contribute (Barreiro-Gomez and Poveda, 2023).

We refer to this new logit dynamic as the generalized pair-wise logit (GPL) dynamic. The existence of a unique solution to this dynamic is proven under a Lipschitz continuity condition of the utility. The $\kappa$-exponential function, particularly its continuity with respect to the shape parameter $\kappa$, yields the parameter continuity of the solution to the dynamic, showing that GPL dynamic continuously reduces to a classical one. We also present a finite difference method to discretize the dynamic (Keliger et al, 2022; Lahkar et al., 2023; Mendoza-Palacios and Hernández-Lerma, 2020).

By invoking the inverse control approach to design the objective function from the control (Komaee, 2021), we heuristically obtain the objective function to be employed in the MFG under the assumption that the transition rate of players' actions is the optimal control. In this way, we can explicitly connect the logit dynamic with the MFG. Furthermore, we obtain the objective function of the MFG leading to the GPL dynamic. We characterize the logit function as the optimal control, a random field in our case, in the MFG. Then, we show that there is a situation where the GPL dynamic is uniquely solvable, while the MFG counterpart may not, when the utility does not satisfy monotonicity (Cecchin et al., 2019; Gomes et al., 2013). This monotonicity is not necessary in the GPL dynamic, which can therefore be formulated under a wider condition. We finally apply the two approaches to fishery management problems in Japan. Consequently, we contribute to the theory and analysis of an evolutionary game model and its MFG counterpart, their computation, and application.

The rest of this paper is organized as follows. **Section 2** presents the GPL dynamic and its analysis. **Section 3** considers an MFG related to the GPL dynamic. **Section 4** is devoted to an application of the two approaches to resource management problems. **Section 5** concludes this study and presents its future perspectives. **Appendix A** has technical proofs of **Propositions** in the main text, and **Appendices B** and **C** have auxiliary survey and computational results.

## 2. Mathematical model

### 2.1 Preliminary setting

The preliminaries here follow the literature (e.g., Lahkar and Riedel, 2015). Throughout this paper, time is denoted as $t \geq 0$. We consider a continuum of actions in a compact set $\Omega = [0,1]$ without loss of significant generality. The $\sigma$-algebra on $\Omega$ is denoted by $\mathcal{B}$. The space of all finite signed measures in $(\mathcal{B},\Omega)$ is denoted by $\mathcal{M}(\mathcal{B},\Omega)$, which associates the total variation norm

$$\|\mu\| = \sup_{g \in G} \left| \int_0^1 g(x)\mu(\mathrm{d}x) \right| \quad \text{for any} \quad \mu \in \mathcal{M}(\mathcal{B},\Omega). \tag{1}$$

Here, $G$ is the collection of measurable functions $g : \Omega \to \mathbb{R}$ with $\sup_{x \in \Omega} |g(x)| \leq 1$. The set of all probability measures on $(\mathcal{B},\Omega)$ is $\mathcal{P}(\mathcal{B},\Omega) = \{\mu \in \mathcal{M}(\mathcal{B},\Omega) |, \ \mu(\Omega) = 1, \ \mu(A) \geq 0 \text{ for any } A \in \mathcal{B}\}$. The distance between $\mu, \nu \in \mathcal{M}(\mathcal{B},\Omega)$ is measured through

$$\|\mu - \nu\| = \sup_{g \in G} \left| \int_0^1 g(x)(\mu(\mathrm{d}x) - \nu(\mathrm{d}x)) \right|. \tag{2}$$

We have $\|\mu\| = 1$ and $\|\mu - \nu\| \leq 2$ for $\mu, \nu \in \mathcal{P}(\mathcal{B},\Omega)$. We define a subset of $\mathcal{M}(\mathcal{B},\Omega)$ as follows: $\mathcal{M}_2(\mathcal{B},\Omega) = \{\mu \in \mathcal{M}(\mathcal{B},\Omega) | \|\mu\| \leq 2\}$. We will use the shorthand notations $\mathcal{M} = \mathcal{M}(\mathcal{B},\Omega)$, $\mathcal{M}_2 = \mathcal{M}_2(\mathcal{B},\Omega)$, and $\mathcal{P} = \mathcal{P}(\mathcal{B},\Omega)$. We have $\mathcal{P} \subset \mathcal{M}_2 \subset \mathcal{M}$.

### 2.2 Utility and logit function

Actions of players are parameterized by $x \in \Omega$. If the distribution of actions is $\mu \in \mathcal{P} \subset \mathcal{M}$, each player choosing an action $x$ gains a utility $U : \Omega \times \mathcal{M} \to \mathbb{R}$. For any $x \in \Omega$ and $\mu, \nu \in \mathcal{M}_2$, $U$ is bounded and satisfies the Lipschitz continuity

$$|U(x,\mu) - U(y,\nu)| \leq L_U (|x - y| + \|\mu - \nu\|). \tag{3}$$

The Lipschitz continuity has been widely assumed in evolutionary game theory (e.g., Perkins and Leslie,

2014; Lahkar and Riedel, 2015; Mendoza-Palacios and Hernández-Lerma, 2022). The utility has been defined on $\Omega \times \mathcal{M}$ but not on $\Omega \times \mathcal{P}$ for a technical reason (see **Proof of Proposition 1** in **Appendix A**).

***Remark 1:*** We focus on the utility $U$ independent from time $t$, whereas the findings presented as follows apply to the cases with time-dependent $U$ if the dependence is sufficiently smooth (i.e., Lipschitz continuous in $[0,+\infty)$; see **Section 4**).

Fix a pair of actions $x, y \in \Omega$. The classical logit function $a : \Omega \times \Omega \times \mathcal{M} \to \mathbb{R}$ represents the rate that the player choosing the action $x$ changes it to $y$ after a small amount of time (e.g., Aiba, 2015; Sandholm and Staudigl, 2018; Leon et al., 2023):

$$a(x, y, \mu) = \frac{1}{1 + \exp\left(\frac{U(x, \mu) - U(y, \mu)}{\eta}\right)} \text{ for any } x, y \in \Omega \text{ and } \mu \in \mathcal{M}, \tag{4}$$

where $\eta > 0$ is the noise intensity; a larger $\eta$ yields a smaller variation of $a(x, y, \mu)$ for $x, y \in \Omega$ owing to blurred information. With this formulation, the rate of jumps from $x$ to $y$ is higher when $U(y, \mu)$ is larger than $U(x, \mu)$. The logit function (4) thus models the utility-driven action choice. A player currently choosing a large utility therefore tends to remain using it, according to the logit function.

We introduce the $\kappa$-exponential function as a hyperbolic regularization of the exponential one (e.g., Kaniadakis et al., 2020):

$$e_\kappa(z) = \begin{cases} \left(\kappa z + \sqrt{\kappa^2 z^2 + 1}\right)^{\frac{1}{\kappa}} & (0 < \kappa \le 1) \\ \exp(z) & (\kappa = 0) \end{cases} \text{ for any } z \in \mathbb{R}. \tag{5}$$

This $e_\kappa$ is a continuous, strictly increasing, and strictly convex function ranging from 0 to $+\infty$. It also satisfies $e_\kappa(z) e_\kappa(-z) = 1$ for any $z \in \mathbb{R}$. The inverse $l_\kappa$ of $e_\kappa$ is

$$l_\kappa(z) = \begin{cases} \dfrac{z^\kappa - z^{-\kappa}}{2\kappa} & (0 < \kappa \le 1) \\ \ln(z) & (\kappa = 0) \end{cases} \text{ for any } z > 0. \tag{6}$$

**Figure 1** plots $e_\kappa(z)$, showing that its gradient increases as $\kappa$ decreases.

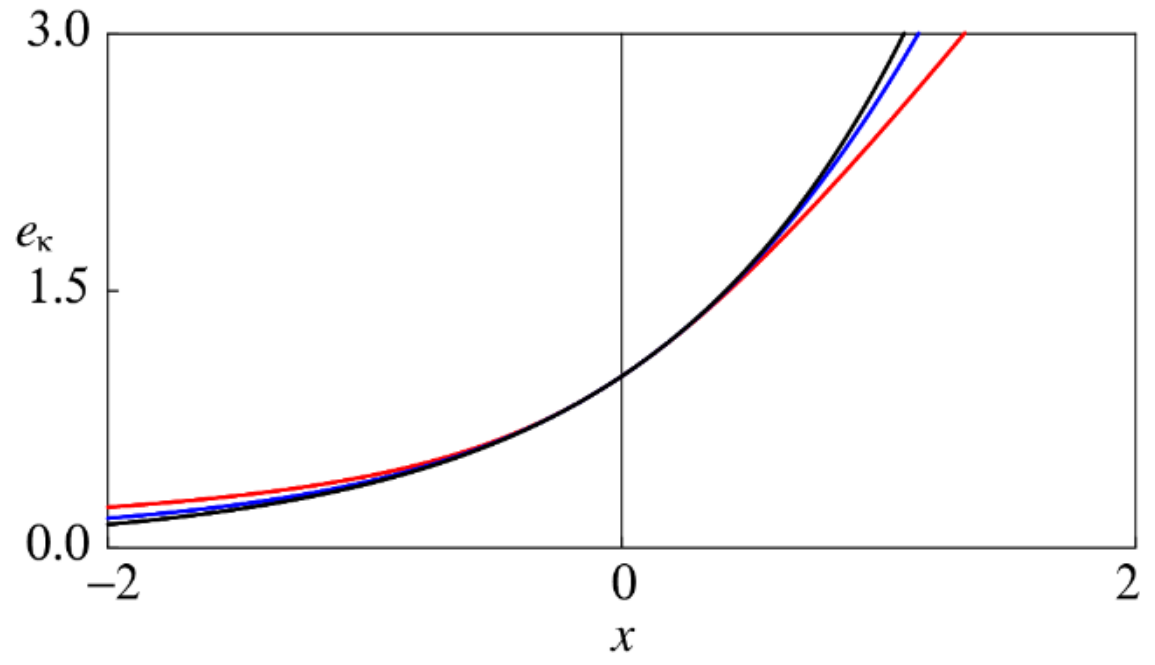

**Figure 1**. Profiles of $e_\kappa(z)$ for $\kappa = 0$ (black), 0.5 (blue), and 1 (red).

Based on the $\kappa$-exponential function, we generalize (4) as

$$a(x, y, \mu) = \frac{1}{1 + e_\kappa\left(\frac{U(x, \mu) - U(y, \mu)}{\eta}\right)} \text{ for any } x, y \in \Omega \text{ and } \mu \in \mathcal{M}. \tag{7}$$

This generalization using the $\kappa$-exponential function implies that the uncertainty that players are faced with has a heavier tail than that assumed in (4). In the sequel, we always assume the logit function $a$ of (7). More specifically, the classical logit function corresponds to a Gumbel as a noise distribution, while that using a $q$-exponential function corresponds to a generalized extreme value distribution having a heavier tail (Nakayama, 2013). Both the $\kappa$- and $q$-exponential functions are hyperbolic generalizations of the exponential one and behave asymptotically as $O\left(y^{1/\kappa}\right)$ and $O\left(y^{1/(q-1)}\right)$, respectively, for a large $y > 0$. They share the same asymptotic growth speed if $\kappa = q - 1$ and $q > 1$. The advantage of using the $\kappa$-exponential function is that it is defined over $\mathbb{R}$; by contrast, the $q$-exponential is defined only for a subset of $\mathbb{R}$ (e.g., Yoshioka, 2024). From this perspective, the $\kappa$-exponential is more useful.

We present the following technical results that will be used later.

***Lemma 1***

*There exists a constant* $L_a > 0$ *depending on* $\kappa \in [0,1]$, $\eta > 0$, *and* $U$ *satisfying (3), such that*

$$\left|a(x_1, y_1, \mu_1) - a(x_2, y_2, \mu_2)\right| \leq L_a\left(|x_1 - x_2| + |y_1 - y_2| + \|\mu_1 - \mu_2\|\right) \tag{8}$$

*for any* $x_1, x_2, y_1, y_2 \in \Omega$ *and* $\mu_1, \mu_2 \in \mathcal{M}_2$.

We also use a technical finding to study continuity of solutions to the GPL dynamic for $\kappa \in [0,1]$.

***Lemma 2***

*For each constant* $M > 0$ *and each* $z \in [-M, M]$, *it follows that* $\frac{\partial e_\kappa(z)}{\partial \kappa}$ *is continuous for* $\kappa \in [0,1]$,

*and further that*

$$\left|\frac{\partial e_{\kappa}(z)}{\partial \kappa}\right| \le C_M \text{ for any } \kappa \in [0,1], \text{ and } \left.\frac{\partial e_{\kappa}(z)}{\partial \kappa}\right|_{\kappa=0} = 0, \tag{9}$$

*where* $C_M > 0$ *is a constant depending only on* $M$.

### 2.3 Generalized pair-wise logit (GPL) dynamic

For any $A \in \mathcal{B}$, the GPL dynamic is given by the following evolution equation that governs a time-dependent probability measure $\mu = \mu(t,\cdot) \in \mathcal{P}$:

$$\frac{\mathrm{d}\mu(t,A)}{\mathrm{d}t} = \int_A \left(\int_\Omega a(z,x,\mu)\mu(t,\mathrm{d}z)\right)\mathrm{d}x - \int_A \left(\int_\Omega a(x,y,\mu)\mathrm{d}y\right)\mu(t,\mathrm{d}x) \text{ for any } t > 0, \tag{10}$$

subject to a prescribed initial condition $\mu(0,\cdot) = \mu_0 \in \mathcal{P}$, where $t$ is the non-dimensional time. The time derivative of $\mu$ is understood in the strong sense with the total variation norm (e.g., Mendoza-Palacios and Hernández-Lerma (2015)). The GPL dynamic (10) describes the evolution of $\mu$ in a way that the local temporal change of $\mu$ (left-hand side) equals the inflow from the other actions (first term in the right-hand side) minus the outflow to the other actions (second term in the right-hand side).

If $\mu$ has a probability density function (PDF), i.e., $\mu(t,\mathrm{d}x) = p(t,x)\mathrm{d}x$, where $p:[0,+\infty)\times\Omega \to [0,+\infty)$ and $\int_\Omega p(t,x)\mathrm{d}x = 1$, then (10) becomes

$$\frac{\mathrm{d}p(t,x)}{\mathrm{d}t} = \int_\Omega a(z,x,\mu)p(t,z)\mathrm{d}z - \left(\int_\Omega a(x,y,\mu)\mathrm{d}y\right)p(t,x) \text{ for any } x \in \Omega \text{ and } t > 0. \tag{11}$$

From this perspective, the connection of the GPL dynamic to jump-driven Markov processes is inferred, where the coefficient $a$ plays the role of a state-dependent jump intensity of a pure-jump Markov process (e.g., Dupret and Hainaut, 2023; Yoshioka et al., 2023; Zhang et al., 2020) and (11) as the associated FP equation. As we will show later, the space-discretized version of the GPL dynamic corresponds to an FP equation of a nonlinear Markov chain. This perspective is a key to the formulation of an MFG version of the GPL dynamic.

***Remark 2:*** The GPL dynamic turns out to be an evolutionary system to maximize the utility $U$ subject to a control cost of the state transition $a$ (see **Section 3**). Particularly, the profile of $a$ is explicitly obtained from the control cost. The GPL dynamic is therefore seen as a costly optimization protocol.

### 2.4 Well-posedness

We prove the unique solvability of GPL dynamic (10), which is our first main result.

***Proposition 1***

*(a) The GPL dynamic (10) subject to a prescribed initial condition* $\mu_0 = \mu(0,\cdot) \in \mathcal{P}$ *admits a unique*

*solution* $\mu = \mu(t,\cdot) \in \mathcal{P}$ *for* $t \in [0,+\infty)$.

*(b) Moreover, if* $\mu_i$ *is the unique solution to the dynamic (10) subject to a prescribed initial condition* $\mu_{0,i} \in \mathcal{P}$ *(*$i = 1,2$*), then there exists a constant* $C > 0$ *independent from* $\mu_{0,1}, \mu_{0,2}$ *such that*

$$\left\|\mu_1(t,\cdot) - \mu_2(t,\cdot)\right\| \leq \left\|\mu_{1,0} - \mu_{2,0}\right\| \exp(Ct), \quad t > 0. \tag{12}$$

We conclude this section with a result on parameter continuity for $\kappa \in [0,1]$, showing that the dynamic (10) is a continuous generalization of the classical one ($\kappa = 0$).

***Proposition 2***

*The unique solution to the GPL dynamic (10) with* $\kappa_i \in [0,1]$*, subject to a common initial condition* $\mu_0$*, is denoted as* $\mu^{(i)}$ *(*$i = 1,2$*). Then, for each* $T > 0$*, it follows that*

$$\lim_{\kappa_1 \to \kappa_2} \sup_{0 \leq t \leq T} \left\|\mu^{(1)}(t,\cdot) - \mu^{(2)}(t,\cdot)\right\| = 0. \tag{13}$$

### 2.5 Discretization

We present a space-discretized version of the GPL dynamic (10). The domain $\Omega$ is divided into $N \in \mathbb{N}$ ($N \geq 2$) cells, where $\Omega_i = \left[(i-1)/N, i/N\right)$ ($i = 1,2,\ldots,N-1$) and $\Omega_N = \left[1-1/N, 1\right]$. The midpoint of $\Omega_i$ is $x_{i-1/2} = (i-1/2)/N$. The dynamic (10) with $A = \Omega_i$ is discretized as

$$\begin{aligned}&\frac{\mathrm{d}\hat{\mu}(t,\Omega_i)}{\mathrm{d}t}\\ &= \frac{1}{N}\sum_{j=1}^{N} a\left(x_{j-1/2}, x_{i-1/2}, \hat{\mu}\right)\hat{\mu}(t,\Omega_j) - \left(\frac{1}{N}\sum_{j=1}^{N} a\left(x_{i-1/2}, x_{j-1/2}, \hat{\mu}\right)\right)\hat{\mu}(t,\Omega_i)\end{aligned} \quad \text{for } i = 1,2,\ldots,N \text{ and } t > 0, \tag{14}$$

where $\hat{\mu}(t,\Omega_i)$ is the discretization of $\mu(t,\Omega_i)$. We used $\int_{\Omega_i} \mathrm{d}x = N^{-1}$ ($i = 1,2,\ldots,N$). The initial condition $\mu(0,\cdot) = \mu_0 \in \mathcal{P}$ is discretized on each $\Omega_i$. The semi-discretized dynamic (14) serves as an approximation to the original one (10), particularly it can be seen as a nonlinear Markov chain with the transition matrix $\left\{a\left(x_{i-1/2}, x_{j-1/2}, \hat{\mu}\right)\right\}_{1 \leq i,j \leq N}$.

To continue, we assume that the utility $U$ has the form

$$U(x,\mu) = g\left(\int_{\Omega} f(x,y)\mu(\mathrm{d}y)\right) \quad \text{for any } x \in \Omega \text{ and } \mu \in \mathcal{P}, \tag{15}$$

where $f: \Omega \times \Omega \to \mathbb{R}$ is continuous on $\Omega \times \Omega$, and $g: \mathbb{R} \to \mathbb{R}$ is Lipschitz continuous such that $\left|g(z_1) - g(z_2)\right| \leq L_g \left|z_1 - z_2\right|$ for any $z_1, z_2 \in \mathbb{R}$, with a constant $L_g > 0$. This utility satisfies the condition (3). If $g(z) = z$, then (15) is the average payoff as conventionally considered in evolutionary games (Sandholm, 2020). Our setting generalizes it in a nonlinear way. We discretize (15) as

$$U\left(x_{j-1/2},\hat{\mu}\right)=g\left(\sum_{k=1}^{N} f\left(x_{j-1/2},x_{k-1/2}\right)\hat{\mu}\left(t,\Omega_k\right)\right),\quad j=1,2,\ldots,N\ . \tag{16}$$

We estimate the difference between the solutions to the original and discretized dynamics. For this purpose, we reinterpret the latter as follows. A piecewise constant approximation of $U$ for any $(x,y)\in\Omega\times\Omega$ and $\mu\in\mathcal{P}$ is set as

$$U_N\left(x,\mu\right)=\sum_{i=1}^{N}U\left(x_{i-1/2},\mu\right)\mathbb{I}_{\Omega_i}\left(x\right)\ \text{ with }\ U\left(x_{i-1/2},\mu\right)=g\left(\sum_{k=1}^{N} f_N\left(x_{i-1/2},x_{k-1/2}\right)\mu\left(\Omega_k\right)\right), \tag{17}$$

($\mathbb{I}_{\Omega_i}\left(x\right)$ is an indicator function that equals 1 if $x\in\Omega_i$, and 0 otherwise) with $f=f_N$ given by

$$f_N\left(x,y\right)=\sum_{i,j=1}^{N} f\left(x_{i-1/2},x_{j-1/2}\right)\mathbb{I}_{\Omega_i}\left(x\right)\mathbb{I}_{\Omega_j}\left(y\right). \tag{18}$$

For any $x,y\in\Omega$, $\mu\in\mathcal{P}$, and $A\in\mathcal{B}$, we set

$$a_N\left(x,y,\mu\right)=\frac{1}{1+e_{\kappa}\left(\dfrac{U_N\left(x,\mu\right)-U_N\left(y,\mu\right)}{\eta}\right)} \tag{19}$$

and

$$I_N\left(\mu,A\right)=\int_A\left(\int_\Omega a_N\left(z,x,\mu\right)\mu\left(\mathrm{d}z\right)\right)\mathrm{d}x-\int_A\left(\int_\Omega a_N\left(x,y,\mu\right)\mathrm{d}y\right)\mu\left(\mathrm{d}x\right), \tag{20}$$

with which the semi-discretized logit dynamic (14) is interpreted as

$$\frac{\mathrm{d}\mu\left(t,A\right)}{\mathrm{d}t}=I_N\left(\mu,A\right),\quad t>0\ . \tag{21}$$

***Remark 3:*** One can generalize (15) to

$$U\left(x,\mu\right)=\sum_{m=1}^{M}g^{(m)}\left(\int_\Omega f^{(m)}\left(x,y\right)\mu\left(\mathrm{d}y\right)\right)\ \text{ for any }\ x\in\Omega\ \text{ and }\ \mu\in\mathcal{P}, \tag{22}$$

where $M\in\mathbb{N}$, and each $f^{(m)}$ and $g^{(m)}$ satisfy the same conditions with $f$ and $g$, respectively.

***Remark 4:*** The existence of unique solutions to (14) follows from the literature (e.g., Proposition 4.1 in Magnus (2023)) as it is a system of ODEs whose coefficients are Lipschitz continuous for each $\mu\left(t,\Omega_i\right)$.

We obtain the following convergence result of the semi-discretized dynamic.

***Proposition 3***

*Solutions to (10) and (21) with the same initial conditions are denoted as $\mu$ and $\mu_N$, respectively. Assume that $\lim_{N\to+\infty}\left\|\mu_N\left(0,\cdot\right)-\mu\left(0,\cdot\right)\right\|=0$ for any $A\in\mathfrak{B}$. Then, it follows that*

$$\lim_{N\to+\infty}\left\|\mu_N\left(t,\cdot\right)-\mu\left(t,\cdot\right)\right\|=0\ \ \textit{at each}\ \ t>0\ . \tag{23}$$

***Remark 5:*** According to **Proof of Proposition 2**, the accuracy of the discretization depends on $\max_{u,v\in\Omega}\left|f_N(u,v)-f(u,v)\right|$. If $\left|f(u_1,v_1)-f(u_2,v_2)\right|\le L_f\left(\left|u_1-u_2\right|+\left|v_1-v_2\right|\right)$ with a constant $L_f>0$ for any $u_1,u_2,v_1,v_2\in\Omega$, then

$$\begin{aligned}\max_{u,v\in\Omega}\left|f_N(u,v)-f(u,v)\right| &= \max_{1\le i,j\le N}\max_{u\in\Omega_i,v\in\Omega_j}\left|f_N(u,v)-f(u,v)\right| \\ &= \max_{1\le i,j\le N}\max_{u\in\Omega_i,v\in\Omega_j}\left|f_N\left(x_{i-1/2},x_{j-1/2}\right)-f(u,v)\right| \\ &\le \max_{1\le i,j\le N}\max_{u\in\Omega_i,v\in\Omega_j}L_f\left(\left|x_{i-1/2}-u\right|+\left|v-x_{j-1/2}\right|\right), \\ &\le \max_{1\le i,j\le N}\max_{u\in\Omega_i,v\in\Omega_j}L_f\left(\frac{1}{2N}+\frac{1}{2N}\right) \\ &= \frac{L_f}{N}\end{aligned} \tag{24}$$

implying the first-order convergence $\left\|\mu_N(t,\cdot)-\mu(t,\cdot)\right\|=O\left(N^{-1}\right)$ for a large $N\in\mathbb{N}$.

## 3. Connection to mean field game

### 3.1 Mean field formulation

We show that the GPL dynamic can be obtained as a version of an MFG to better understand the relationship between them. The discussion in this section is heuristic but gives new insights into the two formulations. For controls as a random field, see the literature (Cartea et al., 2019; Kroell et al., 2023). See **Section 3.3** for the applicability and limitation of the proposed methodology, particularly its theoretical and technical aspects.

Consider a continuous-time jump-driven stochastic process $X=\left(X_t\right)_{t\ge 0}$ with a random field $a_t(\cdot):\Omega\to[0,+\infty)$ ($t\ge 0$). We assume and expect that the solution to (25) is found in a strong (path-wise, right-continuous and has left limit) sense for each initial condition:

$$\mathrm{d}X_t=\mathrm{d}P_t,\quad t>0 \tag{25}$$

with an initial condition $X_0=x\in\Omega$, where $P=\left(P_t\right)_{t\ge 0}$ is a compound jump process with a Lévy density $a_t(\cdot)$ (a product of a PDF of jump size and jump intensity) and jump size $z_t-X_{t-}$ at time $t>0$, where $z_t$ is a jump size valued in $\Omega$. Then, it follows that $X_t$ is almost surely contained in $\Omega$. Set an objective function to be maximized by $a$:

$$\begin{aligned}J(t,x,a) &= \mathbb{E}\left[-\int_t^T\frac{1}{a_s}F\left(a_s\right)\mathrm{d}P_s+\int_t^T G\left(X_s,\mu(s,\cdot)\right)\mathrm{d}s\middle| X_t=x\right] \\ &= \mathbb{E}\left[-\int_t^T\int_\Omega\frac{a_s(y)}{a_s(y)}F\left(a_s(y)\right)\mathrm{d}y\mathrm{d}s+\int_t^T G\left(X_s,\mu(s,\cdot)\right)\mathrm{d}s\middle| X_t=x\right],\quad 0\le t\le T \\ &= \mathbb{E}\left[-\int_t^T\int_\Omega F\left(a_s(y)\right)\mathrm{d}y\mathrm{d}s+\int_t^T G\left(X_s,\mu(s,\cdot)\right)\mathrm{d}s\middle| X_t=x\right]\end{aligned} \tag{26}$$

with $F:\Omega\to\mathbb{R}$ and $G:\Omega\times\mathcal{P}\to\mathbb{R}$, and a fixed terminal time $T>0$, where $\mathbb{E}$ is the expectation with respect to the law $\mu$ of $(t,X_t)$. The objective function (26) implies that the decision-maker, a representative player in the MFG, wants to maximize the net utility (the second term) while maintaining the control cost (the first term) small. The meaning of "cost" will be clarified later. The scaling $F(a_s)$ by $a_s^{-1}$ is to simplify the following analysis. We do not add terminal utility in (26) for simplicity.

We consider $a$ as a Markovian control like the stochastic control theory (Øksendal and Sulem, 2019). Set the admissible set $\mathfrak{A}$ of controls as the collection of predictable random fields $a_t(\cdot):\Omega\to[0,+\infty)$ with respect to a natural filtration generated by $P$. This control modulates intensity and size of jumps of $P$. The HJB equation associated with the maximization problem

$$\Phi(t,x)=\sup_{a\in\mathfrak{A}} J(t,x,a) \quad \text{for any} \quad 0<t<T\,,\quad x\in\Omega \tag{27}$$

is

$$-\frac{\partial\Phi(t,x)}{\partial t}-\sup_a\left\{\begin{array}{l}\int_\Omega a(t,x,y)\big(\Phi(t,y)-\Phi(t,x)\big)\mathrm{d}y\\ -\int_\Omega F\big(a(t,x,y)\big)\mathrm{d}y+G\big(x,\mu(t,\cdot)\big)\end{array}\right\}=0 \quad \text{for any} \quad 0<t<T\,,\quad x\in\Omega \tag{28}$$

subject to the terminal condition $\Phi(T,\cdot)=0$. The supremum in (28) is taken with respect to all measurable functions $a:[0,T]\times\Omega\times\Omega\to[0,1]$. This HJB equation should be coupled with the FP equation that governs time-forward evolution of the probability density $\mu$ (e.g., Example 16 of Bect (2010) where $a^*$ is formally seen as a function of time and space):

$$\frac{\mathrm{d}\mu(t,A)}{\mathrm{d}t}=\int_A\left(\int_\Omega a^*(t,z,x,\mu)\mu(t,\mathrm{d}z)\right)\mathrm{d}x-\int_A\left(\int_\Omega a^*(t,x,y,\mu)\mathrm{d}y\right)\mu(t,\mathrm{d}x) \quad \text{for any} \quad A\in\mathcal{B}\,,\quad t>0\,. \tag{29}$$

Here, we infer the Markovian control $a^*=a^*(t,X_t,\cdot,\mu(t,X_t))$ of the form, as in conventional control problems ($G$ is independent from $a$), as

$$a^*(t,x,y,\mu)=\arg\max_a\left\{\int_\Omega a(t,x,y)\big(\Phi(t,y)-\Phi(t,x)\big)\mathrm{d}y-\int_\Omega F\big(a(t,x,y)\big)\mathrm{d}y\right\}, \tag{30}$$

where the dependence on $\mu$ indicated in (30) considers the implicit dependence of $\Phi$ on it. Given an initial condition $\mu(0,\cdot)$, the MFG as a coupled system (28)–(30) can be solved at least formally, where (28) is integrated with respect to time in a backward manner, while (29) is integrated in a forward manner.

We want to solve the maximization problem (30). Assume that $\Phi$ is bounded and sufficiently smooth, and that the maximization is achieved by an interior solution wherein $0<a<1$. Then, we have

$$\Phi(t,y)-\Phi(t,x)=F'(a), \tag{31}$$

where $F'$ is the derivative of $F$. If $F'$ is strictly increasing except at a finite number of points, then the optimal $a$ should satisfy

$$a^*=a^*(t,x,y,\mu)=F'^{(-1)}\big(\Phi(t,y)-\Phi(t,x)\big). \tag{32}$$

Its right-hand side is a function of the difference of the value function $\Phi$ as a "utility" between $x,y\in\Omega$.

From this perspective, considering the logit function (7) and identifying $U$ with $\Phi$, we design $F$ as

$$F'^{(-1)}(z)=\frac{1}{1+e_{\kappa}\left(-z\eta^{-1}\right)} \quad \text{for any } z\in\mathbb{R}. \tag{33}$$

By setting $z=\eta\ln_{\kappa}\left(\frac{1-u}{u}\right)$ ($0<u<1$), we obtain

$$F'^{(-1)}\left(-\eta\ln_{\kappa}\left(\frac{1-u}{u}\right)\right)=u \quad \text{or equivalently} \quad F'(u)=-\eta\ln_{\kappa}\left(\frac{1-u}{u}\right), \tag{34}$$

and hence

$$F(u)=-\eta\int\ln_{\kappa}\left(\frac{1-u}{u}\right)\mathrm{d}u+F_0 \quad \text{for any } 0<v<1 \tag{35}$$

with some $F_0\in\mathbb{R}$. This $F$ is strictly increasing except at $u=1/2$, has a line symmetry $F(u)=F(1-u)$ with respect to $u\in(0,1)$, and is minimized at $u=1/2$. It is reasonable to assume that there is no cost if the representative player chooses a jump intensity independent from $t,x,y,\mu$, which occurs if $u=1/2$. We thus assume $F(1/2)=0$ in the sequel. Consequently, there will be a cost of controlling the Lévy density if a player wants to design a non-constant one deviating from the uniformly random case $1/2$. This implies that playing the game considering the other players is costly, and hence the cost here can be interpreted as the information cost.

The derived $F$ satisfies $\lim_{x\to 0}F'(|x|)=-\infty$ and $\lim_{x\to 0}F'(1-|x|)=+\infty$, and hence $a$ in (32) is an interior solution. Few values of $\kappa$ give an explicit $F$ that is useful in applications: for any $0<u<1$,

$$\eta^{-1}F(u)=(1-u)\ln(1-u)+u\ln u+\ln 2 \quad \text{if } \kappa=0, \tag{36}$$

$$\eta^{-1}F(u)=1-2\sqrt{u(1-u)} \quad \text{if } \kappa=0.5, \tag{37}$$

$$\eta^{-1}F(u)=-\ln\left(2\sqrt{u(1-u)}\right) \quad \text{if } \kappa=1. \tag{38}$$

The profiles of $F$ in (35) are plotted in **Figure 2**.

***Remark 6:*** The MFG can be numerically dealt with using the forward–backward iteration method that applies to long-time problems (Achdou and Capuzzo-Dolcetta, 2010; Carmona et al., 2023), which we use in **Section 4**. The damping procedure in Barker et al. (2022) is used to obtain numerical solutions with a finite number of iterations. The numerical algorithm is the same as in Barker et al. (2022), except that ours is the explicit finite difference method. See **Appendix C** for the discretization of the MFG.

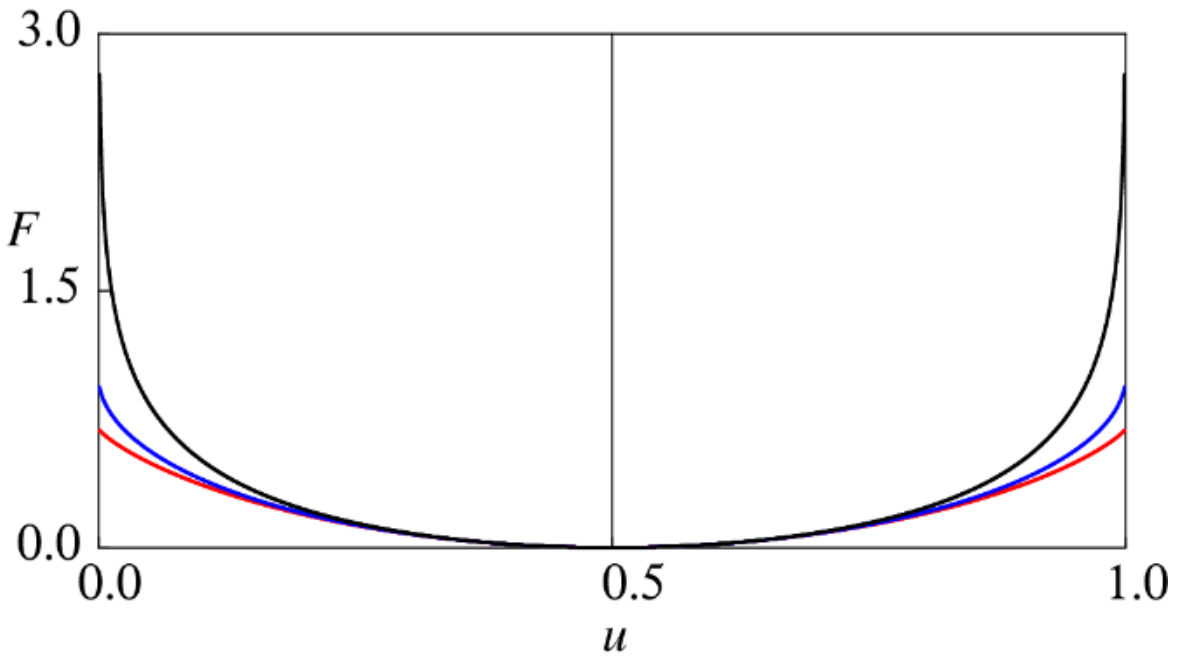

**Figure 2**. Profiles of $F(u)$ for $\kappa$ = 0 (black), 0.5 (blue), and 1 (red). Here, we set $\eta = 1$.

### 3.2 Myopic version

Following the scaling approach (Barker, 2019; Barker et al., 2022), we characterize the GPL dynamic as a version of the MFG. To this aim, for each small increment $h > 0$ of time, we set the scaled objective

$$J_h(t,x,a) = \mathbb{E}\left[-\int_t^{t+h}\int_\Omega F(a_s(y))\,\mathrm{d}y\mathrm{d}s + \frac{1}{h}\int_t^{t+h} G(X_s,\mu(s,\cdot))\,\mathrm{d}s \,\middle|\, X_t = x\right]. \tag{39}$$

The scaling by $h^{-1}$ in the second term, which is motivated from the model predictive control wherein the current control applies, is designed using the prediction for the near future (Barker, 2019). To proceed further, we assume $G = U$ with (15) and that the law of $X$ is stationary and, again, that the control $a$ is Markovian. We infer an approximation of the right-hand side of (39):

$$\begin{aligned}
&\mathbb{E}\left[-\int_t^{t+h}\int_\Omega F(a(X_s,y))\,\mathrm{d}y\mathrm{d}s + \frac{1}{h}\int_t^{t+h} G(X_s,\mu(s,\cdot))\,\mathrm{d}s \,\middle|\, X_t = x\right] \\
&= \mathbb{E}\left[-\int_t^{t+h}\int_\Omega F(a(X_s,y))\,\mathrm{d}y\mathrm{d}s + \frac{1}{h}\int_t^{t+h} U(X_s,\mu(s,\cdot))\,\mathrm{d}s \,\middle|\, X_t = x\right]. \\
&\approx h\mathbb{E}\left[-\int_\Omega F(a(X_t,y))\,\mathrm{d}y + \frac{1}{h}U(X_{t+h},\mu(t,\cdot)) \,\middle|\, X_t = x\right]
\end{aligned} \tag{40}$$

For bounded and sufficiently smooth $\phi = \phi(t,x) : [0,+\infty) \times \Omega \to \mathbb{R}$, Itô's lemma (e.g., Theorem 1.14 of Øksendal and Sulem (2019)) implies

$$\mathbb{E}\left[\phi(X_{t+h}) \middle| X_t = x\right] \approx \phi(x) + h\int_\Omega a(x,y)(\phi(y) - \phi(x))\,\mathrm{d}y\,. \tag{41}$$

By (15), if $U(\cdot,\mu)$ is sufficiently smooth, then we infer

$$\begin{aligned}
&\frac{1}{h}\mathbb{E}\left[-\int_t^{t+h}\int_\Omega F(a(X_s,y))\,\mathrm{d}y\mathrm{d}s + \frac{1}{h}\int_t^{t+h} U(X_s,\mu)\,\mathrm{d}s \,\middle|\, X_t = x\right] \\
&\approx \mathbb{E}\left[-\int_\Omega F(a(x,y))\,\mathrm{d}y + \frac{1}{h}\left(U(x,\mu) + h\int_\Omega a(x,y)(U(y,\mu) - U(x,\mu))\,\mathrm{d}y\right) \,\middle|\, X_t = x\right]. \\
&= \frac{1}{h}U(x,\mu) + \int_\Omega \left\{-F(a(x,y)) + a(x,y)(U(y,\mu) - U(x,\mu))\right\}\mathrm{d}y
\end{aligned} \tag{42}$$

The maximizing $a$ of (39) would be close to that of (42) if $h$ is small. As in the derivation procedure for (35), the maximizing $a=\hat{a}$ of the last integral of (42) is

$$\hat{a}(x,y,\mu)=F'^{(-1)}\left(U(y,\mu)-U(x,\mu)\right), \tag{43}$$

which is identical to (7) if $F$ is given by (35). Consequently, the myopic version of the MFG based on the ansatz mimicking the model predictive control yields the logit function (7), and hence the GPL dynamic.

The approach from the MFG suggests that the logit function (7) can be interpreted as the optimized jump intensity of the process $X$ according to a scaled control problem. The functional form and the derivation procedure, namely (43), suggest that if the current utility is high such that $U(x,\mu)>>U(y,\mu)$, then the logit function is close to 0. Then, the player should not change the current action. We thus found an explicit connection to the logit function and a controlled pure-jump stochastic process; the latter can be understood as a trajectory of actions chosen by a representative player. This kind of explicit connection between a logit dynamic and an MFG has not been reported so far.

### 3.3 Comparison and discussion

We consider that a key in our methodology is specifying a relevant and explicit coefficient $F$ (i.e., (36)-(38)) in the MFG. Fortunately, we could find the explicit $F$ for the selected values of $\kappa=0,\ 0.5,\ 1$ but not for the other values. The situation may get worse for the other logit functions not studied here. The calculation to connect the logit function and the coefficient $F$ is very technical, and its applicability to the other evolutionary game models should be investigated in the future. Another limitation of this study is that presents a linkage between the GPL dynamic and its MFG version only heuristically. To rigorously address this issue, possibly based on stochastic calculus, one may need to start from convergence of a many-player model to infinite-player one as discussed in the literature for both the evolutionary and mean field game models (e.g., Bayraktar et al., 2022; Laurière and Tangpi, 2022). Our contribution does not address this issue; instead, we present some tractable cases where the coefficients of the MFG can be found explicitly. Obtaining a more theoretical connection between an MFG to a GPL dynamic is beyond the scope of this study, which is why we present some computational applications. Deeper theoretical investigations will need an effort from a purely mathematical side where the class of utility and the range of model parameters may have to be more narrowly specified. Despite these limitations, due to dealing with specific problems in detail, we believe that our methodology would be give hints and insights into deeper understanding between evolutionary and mean field game models from moth theoretical and engineering sides.

The GPL dynamic is a simpler model than the MFG version because the former has a smaller number of parameters than the latter. Analysis of stationary states can be more efficiently conducted by the GPL dynamic. By contrast, the MFG is potentially able to deal with more complicated problems where influences of foresights of players are of interest.

Guaranteeing the unique solvability of an MFG problem often requires the monotonicity. Notice that we are dealing with a maximization problem (past research (e.g., Cohen and Zell, 2023; Gomes et al., 2013) dealt with minimization problems):

$$\int_{\Omega}\left(U(x,\mu)-U(x,\nu)\right)\left(\mu(\mathrm{d}x)-\nu(\mathrm{d}x)\right)\leq 0 \quad \text{for any} \quad \mu,\nu\in\mathcal{P} \quad \text{and} \quad x\in\Omega . \tag{44}$$

This condition is satisfied if $U$ is of the potential game type often considered in applications (e.g., Friedman and Ostrov, 2013; Leonardos and Piliouras, 2022; Tan et al., 2023). Phenomenologically, the monotonicity penalizes congestion (Delarue and Tchuendom, 2020), namely, the concentration of actions that may eventually lead to a bifurcation and, hence, non-uniqueness of solutions. This monotonicity is not needed for the unique solvability of the GPL dynamic (**Proposition 1**). Without the monotonicity, the existence follows if the Lipschitz continuity is satisfied and the action space has been discretized as in **Section 2.5**; however, uniqueness would not necessarily apply (Gomes et al., 2013). This difference between the two problems sharing the same set of parameter values (except for $T$) is caused by the different assumptions in their perspectives against the decision-making. In an evolutionary game that includes (generalized pair-wise) logit dynamics, each player chooses actions based only on the current information without any foresights. By contrast, in the MFG, he/she chooses actions to maximize an objective function, which is a conditional expectation and therefore a prediction. Indeed, the model predictive control is more myopic than the original MFG, as inferred from the analysis in **Section 3.2**. Later, we show an example where the GPL dynamic successfully computes a solution, while the corresponding MFG does not numerically converge, possibly owing to the absence of monotonicity.

The GPL dynamic is also related to the myopic adjustment process (Neumann, 2023b) as a reduced model of an MFG with a discounted objective function. As already discussed, a common drawback of MFGs is that each game involves a pair of evolution equations that must be solved in opposite directions in time. This aspect would be an especially serious matter about the feasibility of MFGs in applications because it possibly implies a huge computational time. Focusing on a finite-state system, Neumann (2023b) proposed to reduce its MFG to a nonlinear Markov chain whose coefficients are such that an optimal strategy by the reduced model is a mean field equilibrium, to ensure that they are consistent with each other at a stationary state. Our model reduction of an MFG to a logit dynamic is theoretically different from this approach, albeit both focus on a stationary state.

Finally, we emphasize a relation of the proposed dynamic to graphon games (Parise and Ozdaglar, 2019). A graphon is a continuum limit of adjacency matrix among players. In our context, a graphon can be a function $\vartheta:\Omega\times\Omega\to[0,+\infty)$ such that the connection between the pair of players $(x,y)$ is higher for a larger $\vartheta$. Through the introduction of a graphon, the logit function (7) can be extended as

$$a(x,y,\mu)=\frac{\vartheta(x,y)}{1+e_{\kappa}\left(\frac{U(x,\mu)-U(y,\mu)}{\eta}\right)} \quad \text{for any} \quad x,y\in\Omega \quad \text{and} \quad \mu\in\mathcal{P} . \tag{45}$$

In this way, our GPL dynamic uses a uniform graphon $\vartheta\equiv 1$ (Petrov, 2023). Investigating models with a nonuniform one will be an interesting topic, which will be addressed elsewhere.

As a final remark, a limitation of the GPL dynamic compared to other evolutionary game models is that it cannot reproduce a stationary solution given by a (superposition of) Dirac's delta, which is a pure strategy. Intuitively, this limitation is explained by the fact that under a stationary state we have $p(x) = \left(\int_\Omega a(x, y, \mu) \mathrm{d}y\right) \int_\Omega a(z, x, \mu) p(z) \mathrm{d}z$ , which does not represent any pure strategies. The same drawback applies to the MFG version. By contrast, replicator dynamics can reproduce stationary pure strategies (e.g., Mendoza-Palacios and Hernández-Lerma, 2015).

## 4. Applications

### 4.1 Study target

We conduct studies on the application of the GPL dynamic MFG to two fisheries management problems in inland fisheries in Japan. Both applications focus on the fish *Plecoglossus altivelis altivelis* (*P. altivelis*) in the Hii River, San-in area, Japan. This fish species has been a pivotal inland fishery resource in the country, occupying 9.6% of the total fish catch (Ministry of Agriculture, Forestry and Fisheries, 2023), and has been highlighted for its significant economic value (Inui et al., 2021). The fish is famous among recreational fishing enthusiasts during summer (Hata and Otake, 2019; Yoshioka, 2023; Yoshimura et al., 2021), and its population decrease has recently been a critical concern owing partly to climate changes (Nagayama et al., 2023). Inland fishery cooperatives have therefore been searching for ways to sustain the fish population while supplying fishery resources for recreational fishing and revitalizing the regional economy. The Hii River Fishery Cooperative (HRFC) has been overseeing inland fishery resources in the mid- to up-stream reaches of the Hii River. We have been collaborating with HRFC since 2015 to collect ecological data related to *P. altivelis* in the Hii River.

### 4.2 First application

The first application concerns recreational fishing for one fishing season of *P. altivelis* (July 1 to October 30 in a year), which is a problem in a monthly scale. In this application, $x \in \Omega$ is the arrival intensity to the Hii River with the unit 1/day, and thus each angler must buy fishing tickets depending on $x$. In this problem setting, each angler arrives at the river at most once in a day. We consider a time-dependent utility to be able to account for the fishing utility gained by catching larger fish.

The utility considered in this application is the time-dependent one

$$U_t(x, \mu) = \int_\Omega \left\{ -A_1 (1 + A_2 y) x + \left( \frac{W_t}{W_{\max}} x \right)^{A_3} \right\} \mu(\mathrm{d}y) \quad \text{for any} \quad x \in \Omega, \ t \geq 0, \ \mu \in \mathcal{P}, \tag{46}$$

where $A_1, A_2, A_3 \geq 0$ and the body weight $W_t > 0$ of the fish. The first term in the right-hand side of (46) represents the arrival cost, assuming that each player (angler) pays the daily fee (Hii River Fishery Cooperative, 2023), whereas the coefficient $A_2 y$ represents the added cost due to harvesting excessive to the mean. This added cost is a manifestation of a sustainability concern such that a larger mean arrival rate of anglers leads to a larger net cost of harvesting per angler, which will effectively protect the fishery resource from overharvesting. The second term represents the utility gained by harvesting, which increases with respect to the mean body weight $W_t$ of *P. altivelis* at time $t > 0$. The division by the maximum body weight $W_{\max} > 0$ is to normalize this term. This utility is monotone at each time $t \geq 0$.

The body weight $W_t$ is modeled by the logistic model (Yoshioka and Yaegashi, 2016)

$$\frac{W_t}{W_{\max}} = \left\{1 + \left(\frac{W_{\max}}{W_0} - 1\right)\exp(-rt)\right\}^{-1}, \quad t \geq 0 \tag{47}$$

with initial body weight $W_0 > 0$ and growth rate $r > 0$. The least-squares fit against the data collected in 2022 gives the parameter values $W_{\max} = 106$ (g), $W_0 = 17.0$ (g), and $r = 0.0223$ (1/day). Here, the initial time 0 is the beginning of May 1 (see **Appendix B**). We set $A_2 = 1$ and $A_1 = A_3 = 1/2$ unless otherwise specified, with which each term is balanced, and the computed transient PDF is almost unimodal.

We compare the GPL dynamic and MFG for the 120 days from the beginning of July 1, covering the major fishing season of *P. altivelis* in the Hii River. The terminal time is $T = 120$ (day), and the time increment for the temporal discretization is 0.01 (day). The time $t$ in this case is not normalized and has the dimension of day, suggesting that players would be able to update their actions in a daily timescale. We examine the three values of $\kappa = 0.0,\ 0.5,\ 1.0$ for which the cost $F$ is explicitly available and a larger $\eta = 0.1$ and smaller $\eta = 0.01$. The MFG is computed using the damping factor $\omega = 0.250$.

**Figures 3–4** show the computed transient PDFs of the GPL dynamic for different values of $\kappa$ with $\eta = 0.1$ and $\eta = 0.01$, respectively. **Figure 3** suggests that for $\eta = 0.1$, the PDF profiles are comparable among the different values of the shape parameter $\kappa$, whereas **Figure 4** shows that they are distinguishable, with different magnitudes for the height of a single peak, specifically decreasing with respect to $\kappa$. This is in accordance with the intuition that $\kappa$ governs the tail of the uncertainty, and a larger $\kappa$ implies a larger variance. **Figures 5** shows the computed PDFs of the MFG with $\eta = 0.01$, demonstrating that the PDF profiles are sharper than those of the GPL dynamic but have a similar parameter dependence on $\kappa$: the peak height decreases with respect to $\kappa$. The computational results suggest that the short-sightedness of the GPL dynamic recommends more variable action profiles to be optimal.

We also investigate influences of $A_2$ on the dynamic to better understand the role of the sustainability concern. **Figure 6** shows the computed transient PDFs of the GPL dynamic for $\eta = 0.01$ and $\kappa = 0.5$, where the parameter values of $A_2$ are varied to analyze the influence of the sustainability concern on the dynamics. **Figure 7** shows the corresponding results of the MFG. In both cases, the absence of the sustainability concern leads to a profile with its maxima placed closer to the boundary $x = 1$ than that with the sustainability concern. The dependence of the optimal actions on the body weight is more visible in this case, wherein anglers arrive at the river more often in the later fishing season. This potentially leads to a smaller number of fish contributing to spawning in the coming autumn. This kind of weight-dependent action profile is less clearly seen when $A_2 = 1$, and is further suppressed for the larger value $A_2 = 2$. The computational results suggest the importance of considering sustainability in the fisheries management because it was found to lead to qualitatively different action profiles, and hence different intensities of fishing pressure from anglers.

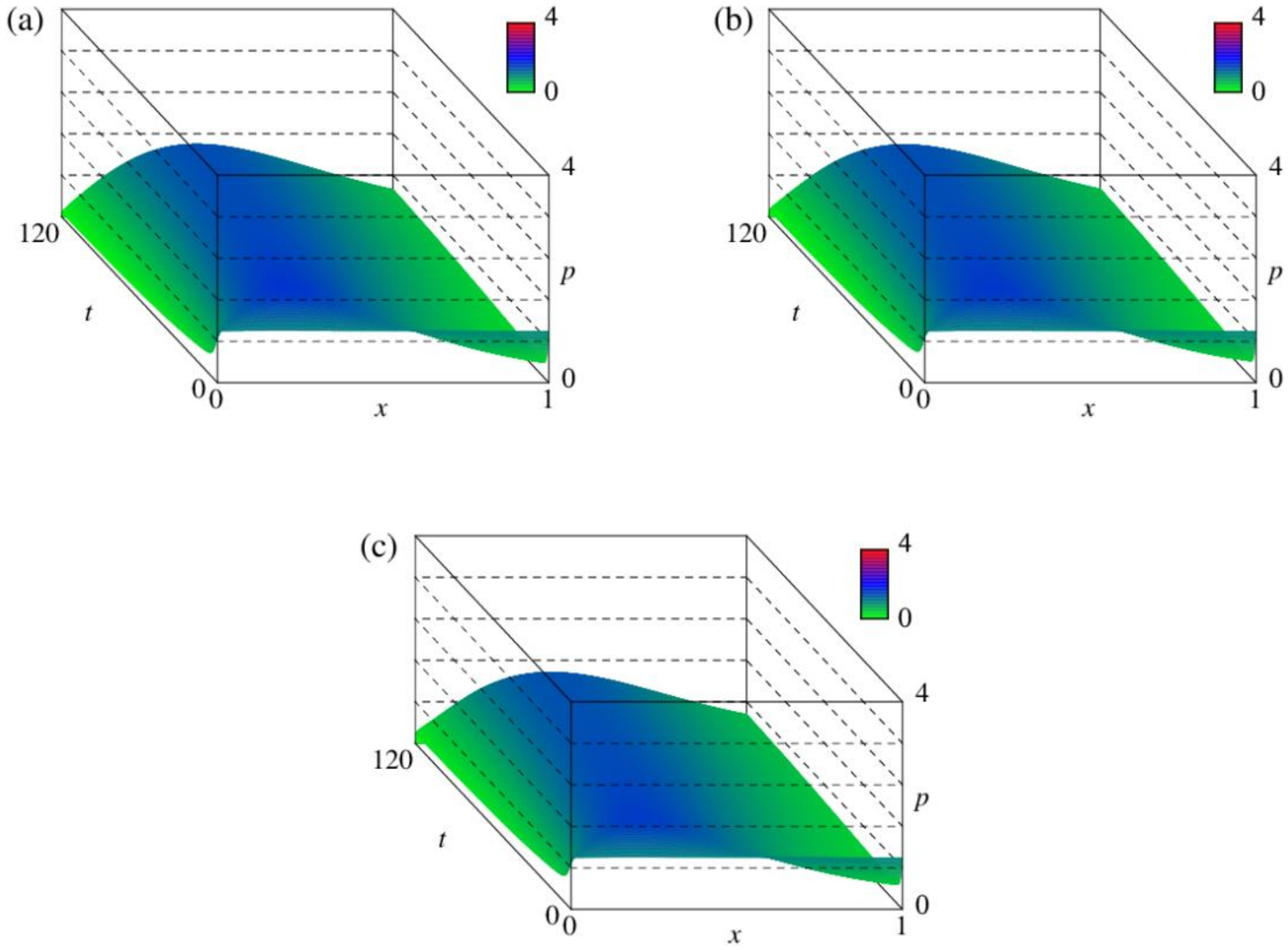

**Figure 3.** Computed PDFs of GPL dynamic with $\eta = 0.1$: (a) $\kappa = 0.0$, (b) $\kappa = 0.5$, (c) $\kappa = 1.0$.

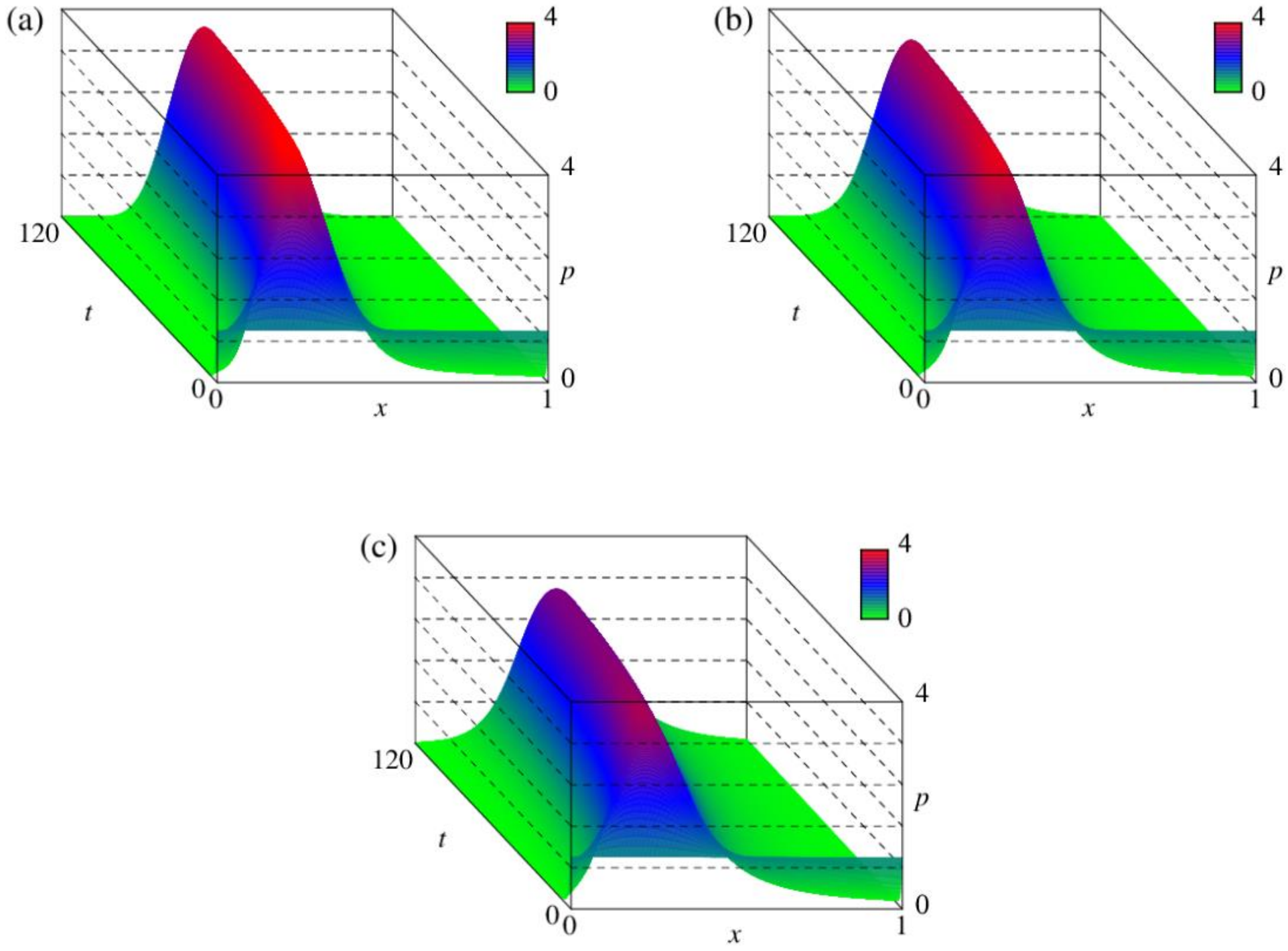

**Figure 4.** Computed PDFs of GPL dynamic with $\eta = 0.01$: (a) $\kappa = 0.0$, (b) $\kappa = 0.5$, (c) $\kappa = 1.0$.

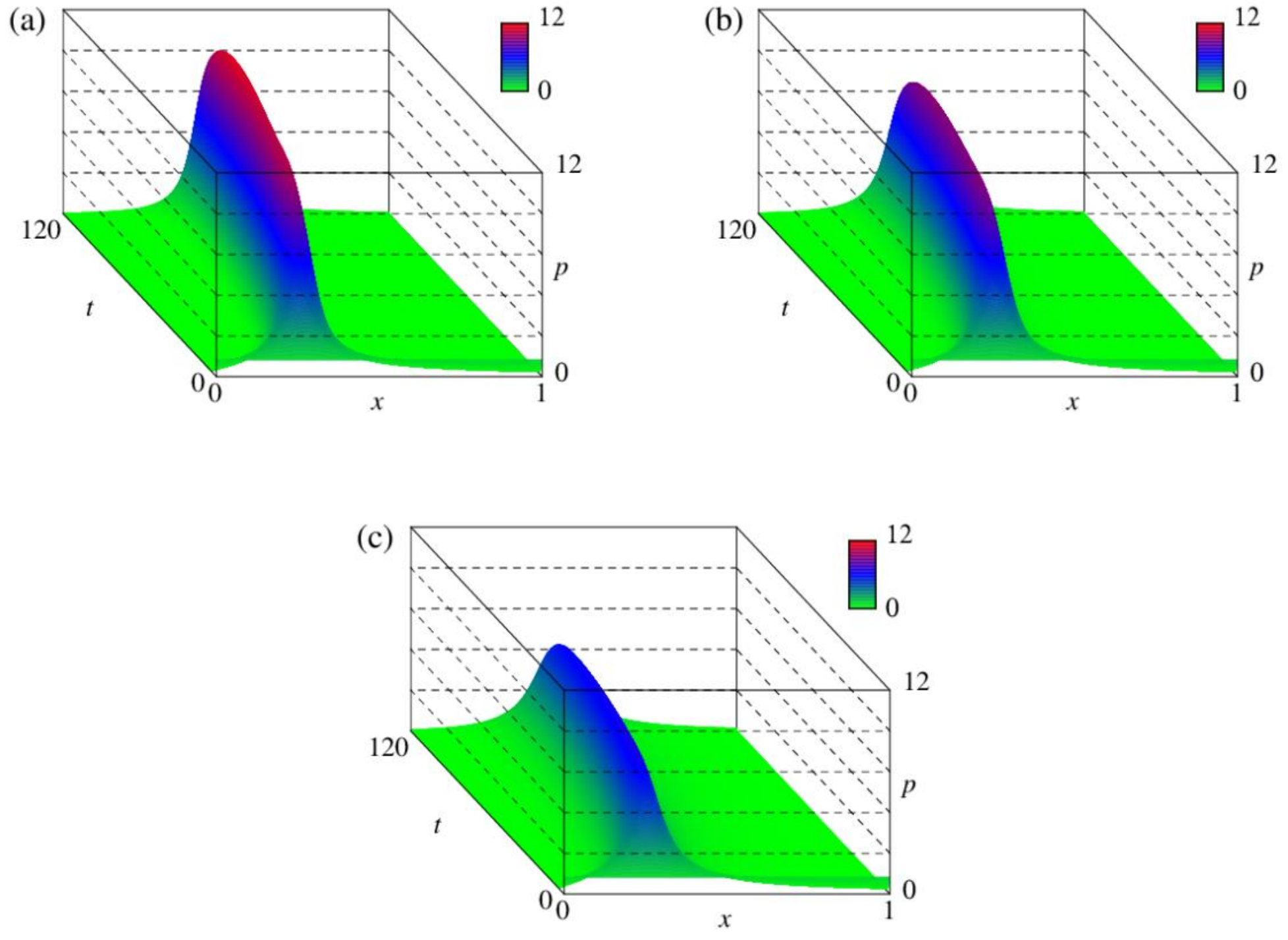

**Figure 5.** Computed PDFs of MFG with $\eta = 0.01$: (a) $\kappa = 0.0$, (b) $\kappa = 0.5$, (c) $\kappa = 1.0$.

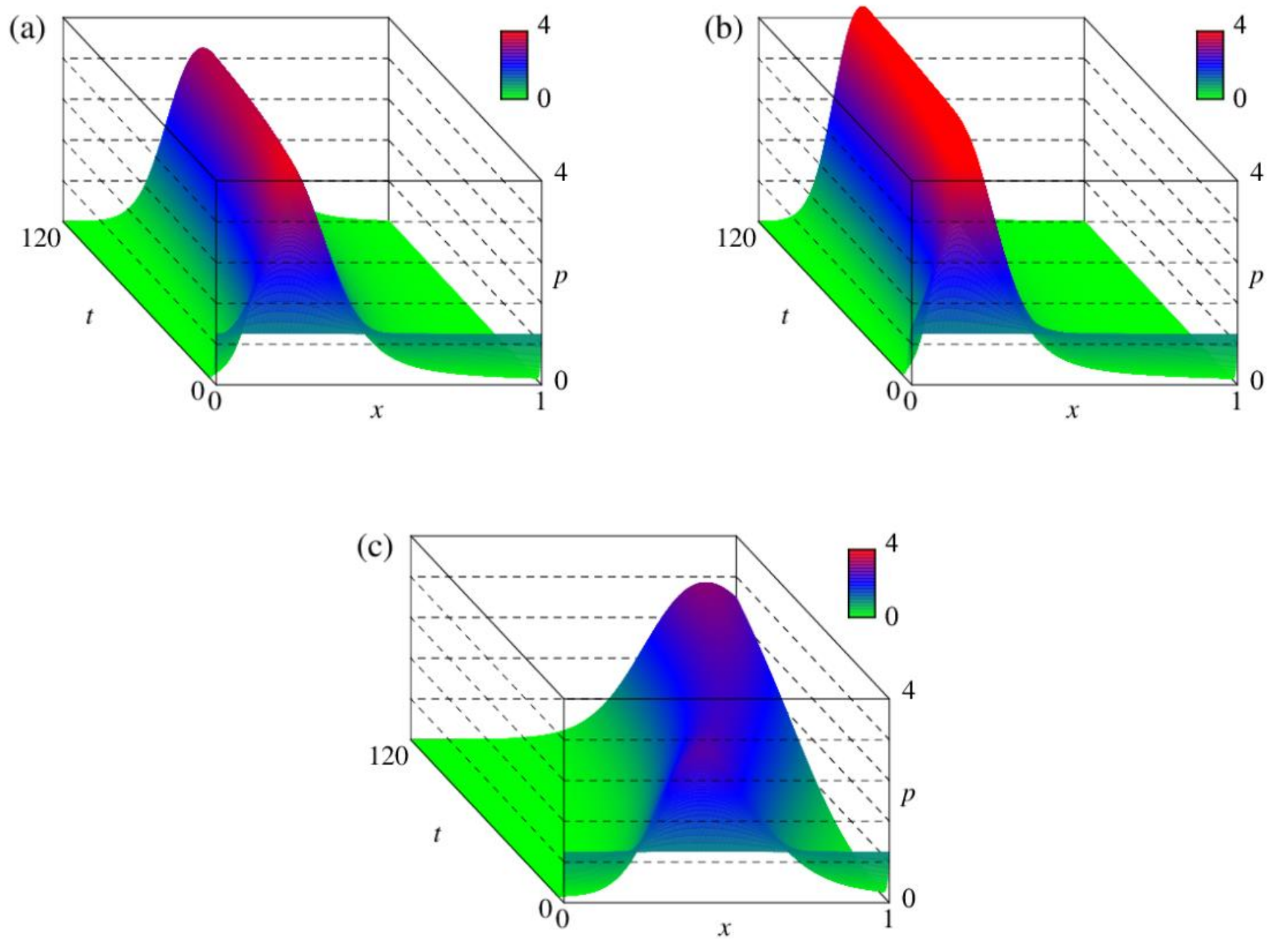

**Figure 6.** Computed transient PDFs of GPL dynamic with $\eta = 0.01$ and $\kappa = 0.5$, (a) nominal case $A_2 = 1$, (b) greater sustainability concern $A_2 = 2$, (c) no sustainability concern $A_2 = 0$.

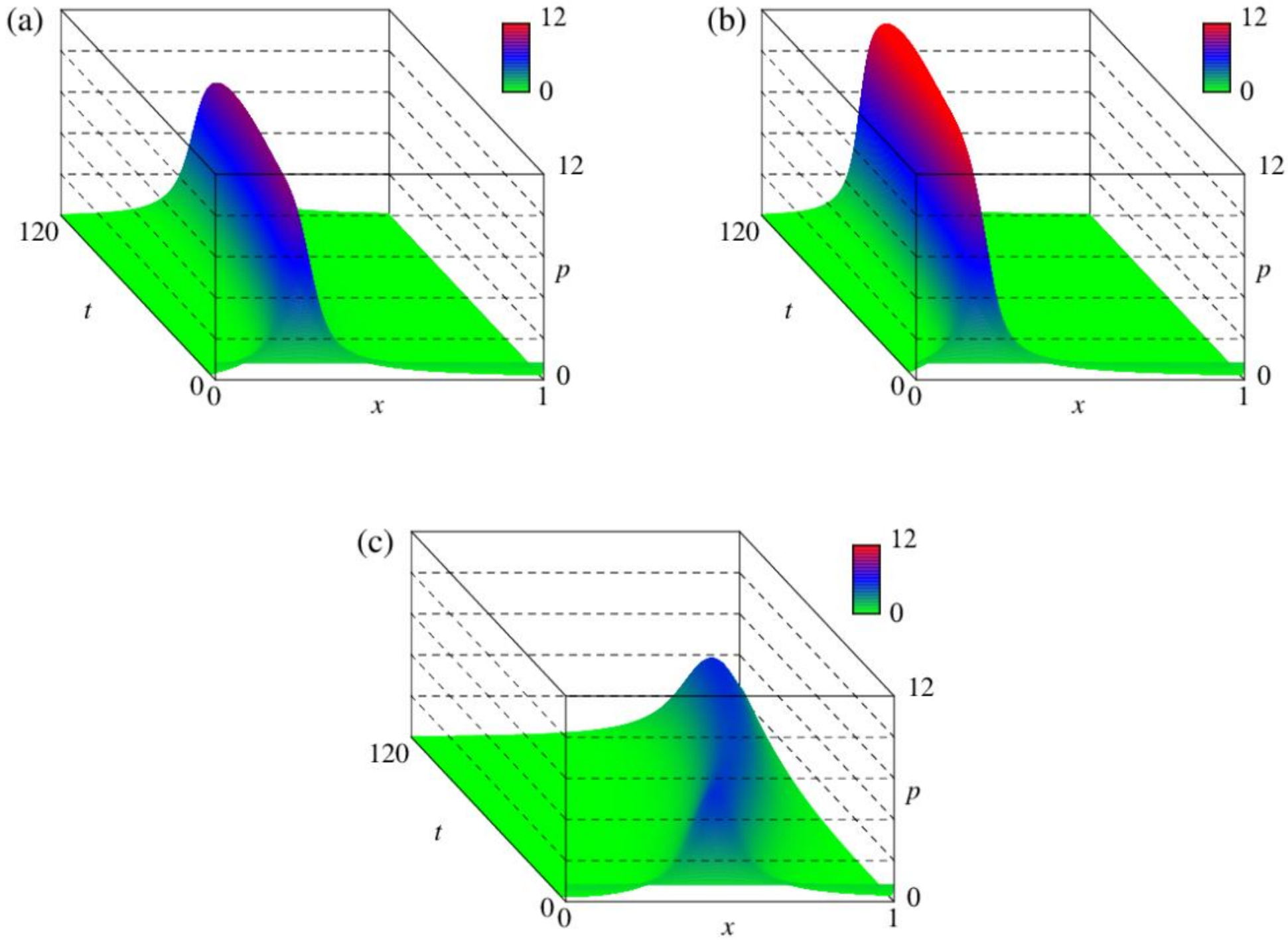

**Figure 7.** Computed transient PDFs of MFG with $\eta = 0.01$ and $\kappa = 0.5$, (a) nominal case $A_2 = 1$, (b) greater sustainability concern $A_2 = 2$, (c) no sustainability concern $A_2 = 0$.

## 4.3 Second application

For the second application, we consider the annual Toami competition. This competition is typically held around the end of July to the beginning of August. Participants are registered in pairs and, on the competition day, compete on the total *P. altivelis* fish catch for two hours (typically 10:00 am to 12:00 am). All participants must use casting nets under a common regulation. The top few pairs (10% to 20% of participants) are awarded by the HRFC.

Compared to the earlier application, the problem considered here is less related to sustainable resource management but has a utility with a unique term being the award scheme employed by the HRFC. Here, we consider $\Omega$ as the normalized total fish catch. The utility considered in this study is that of Yoshioka et al. (2024); for any $x \in \Omega$ and $\mu \in \mathcal{P}$, we set

$$U(x,\mu) = \int_{\Omega} \left\{ -B_1 x^2 + B_2 |x-y|^{B_3} \right\} \mu(\mathrm{d}y) + B_4 \max\left\{ \alpha - \int_x^1 \mu(\mathrm{d}y), 0 \right\}, \quad (48)$$

where $B_1, B_2, B_3, B_4 \geq 0$ and $\alpha \in (0,1)$. The first term of (48) is the harvesting cost, the second term the individual difference to induce a competition, and the third term the awarding scheme by the HRFC; the top $100\alpha\%$ pairs are awarded. This utility was proposed heuristically in Yoshioka et al. (2024) to model the award scheme by the HRFC. We regularize the third term to obtain a Lipschitz continuous utility (Yoshioka et al., 2024):

$$\int_x^1 \mu(\mathrm{d}y) \to \int_0^1 \max\left\{ 0, \min\left\{ 1, \frac{y-x+\varepsilon}{\varepsilon} \right\} \right\} \mu(\mathrm{d}y) \text{ with } 0 < \varepsilon << 1. \quad (49)$$

So far, we did not find the inequality (44) nor its complementary one for the utility (48) and its regularized version, so the monotonicity of this utility is open. Nevertheless, due to the existence of the third term, this utility is worth investigating because it has the unique form reflecting the actual award scheme.

Without loss of generality, we set $B_4 = 1$. Considering the award scheme that the HRFC has used, we set $\alpha = 0.2$. Further, we set $B_3 = 1$ to simplify the problem. We fit the following parameters to the empirical data assuming stationarity: $B_1 = 0.1i$ and $B_2 = 0.1j$ ($1 \leq i, j, \leq 30$), $\eta = 0.01,\ 0.05,\ 0.1$, and $\kappa = 0.1,\ 0.5,\ 1.0$. We then find the best triplet $(B_1, B_2, \kappa)$ for each $\eta$ such that the following error metric based on a moment matching is minimized: $(\mathrm{Ave_e} - \mathrm{Ave_m})^2 + (\mathrm{Std_e} - \mathrm{Std_m})^2$, where Ave=Average, Std=Standard deviation, "e" represents the empirical value, and "m" the modeled value. Stationary solutions to the GPL dynamic are computed using semi-discretization with $N = 500$ cells, which is discretized in time using the forward Euler method with the time increment 0.1. The damping factor used in computing the MFG is $\omega = 0125$. We use $\varepsilon = 1/N$ in (49).

The parameter set identified is summarized in **Table 1**. **Figure 8** compares the empirical and theoretical stationary PDFs to the GPL dynamic using the parameter values in **Table 1**. The performance of the identified models improves as the noise intensity $\eta$ increases, and all of them predict the bimodal PDFs. The identified sizes of the parameters $B_1, B_2$ increase as $\eta$ increases, suggesting that the last term of (48) becomes less important for a model with a larger $\eta$. The disappearance of the right maxima as $\eta$ increases suggests that the players become less sensitive to cost/benefit by harvesting a larger number of fish under the uncertain and, hence, less informative environment. We also compute the MFG counterpart using the computational resolution $\Delta t = 0.01$. The computation has been performed for the three parameter sets in **Table 1** and terminal time $T$ = 40, 80, and 160. We did not obtain convergence of the numerical solutions with the smallest $\eta = 0.01$ against all three values of terminal time. A smaller $\eta$ implies a stronger congestion (i.e., more concentrated maxima of $p$), possibly resulting in a more unstable computational procedure with larger solution gradients.

For $\eta = 0.05$ and $\eta = 0.1$, we obtained converged numerical solutions as shown in **Figures 9–10**. Convergence analyses of the GPL dynamic and MFG are presented in **Appendix C**. For $\eta = 0.05$, **Figure 11** compares the stationary PDF for the GPL dynamic and that with the corresponding MFG solutions at time $t = T/2$ for $T$ = 40, 80, and 160. **Figure 12** shows the results for $\eta = 0.1$. For the MFG, the PDF is considered to be close to a "quasi-stationary state" at time $t = T/2$, around which it is almost invariant with respect to time, and the value function $\Phi$ has a common constant increment in $\Omega$ (Achdou and Capuzzo-Dolcetta, 2010). The time instance $t = T/2$ is therefore suitable for comparing the GPL dynamic and MFG (see **Appendix C**). For both $\eta = 0.05$ and $\eta = 0.1$, the computational results suggest that for the MFG, a more conservative strategy with more concentrating actions on a smaller $x$ is optimal. The short-sightedness of the GPL dynamic therefore recommends an excessive cost in conjunction with high return action for the Toami competition. Contrary to the findings of Barker et al. (2022), the PDF of the GPL dynamic is more diffused than that of the MFG. However, directly comparing their and our models is difficult because of fundamental differences owing to assuming different noise processes, different objective functions and, more importantly, different evolutionally game models. Systematically comparing the evolutionary game and the associated MFG for generic stochastic system dynamics would require the development of a unified mathematical theory along with further case studies.

Finally, we discuss the convergence of the iterative procedure for computing the MFG. **Figure 13** shows the convergence histories for the MFG, showing that the convergence is achieved within a few hundred steps, and that an exponential convergence is found for each computational case. In **Figure 13**, the case $\eta = 0.05$ and $T = 160$, which has the smallest noise intensity and largest terminal time, exhibits a slower convergence than those observed for the other cases. Moreover, its convergence history involves oscillations between 50 and 100 iterations, after which the convergence speed slows down. We found that specifying a sufficiently large $T$ or a sufficiently small $\eta$ results in convergence failure in computing the MFG due to oscillations of numerical solutions between successive iterations. On the other hand, the case of $\eta = 0.05$ and $T = 160$ would be the intermediate case, where the convergence of numerical solutions can be achieved in a somewhat unstable manner.

**Table 1.** Identified parameter set for each $\eta$. Number in "( )" is the relative error with respect to the empirical value. Empirical Ave and Std are 0.32471 and 0.30352, respectively.

| $\eta$ | $B_1$ | $B_2$ | $\kappa$ | Ave | Std |
|---|---|---|---|---|---|
| 0.01 | 0.2 | 0.1 | 1.0 | 0.30242 (2.23E-02) | 0.31250 (8.98E-03) |
| 0.05 | 0.4 | 0.4 | 0.0 | 0.32360 (3.91E-03) | 0.29960 (1.11E-03) |
| 0.10 | 1.3 | 2.3 | 1.0 | 0.32386 (3.13E-05) | 0.30355 (8.47E-04) |

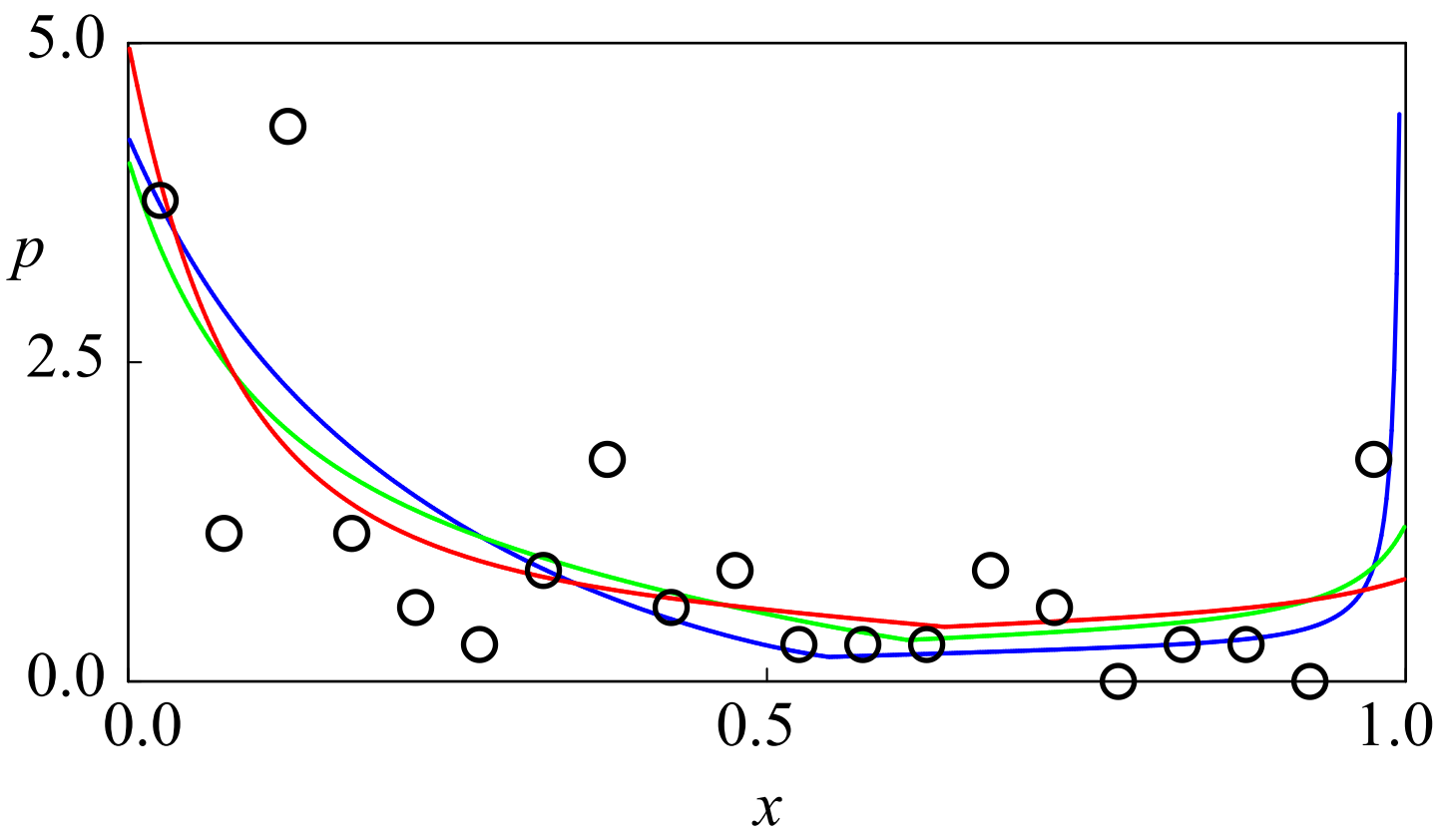

**Figure 8.** Empirical (circles) and theoretical stationary PDFs (curves) for GPL dynamic. Theoretical results based on $\eta$ = 0.01 (blue), 0.005 (green), 0.1 (red).

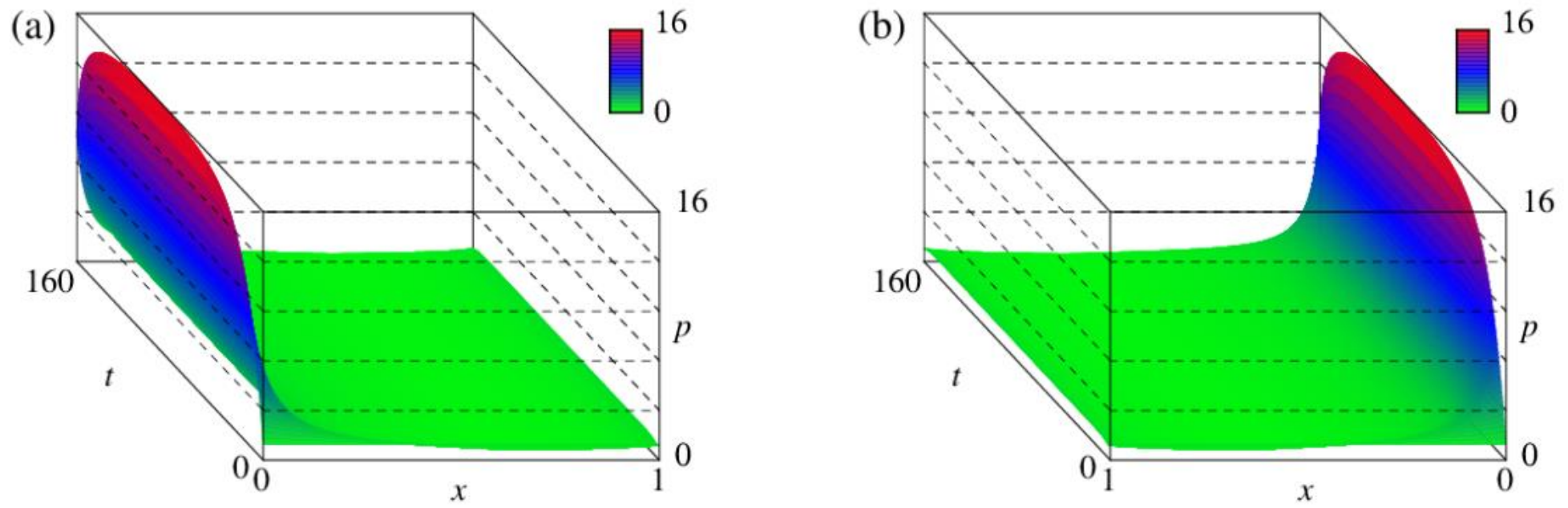

**Figure 9.** Computed PDF of GPL dynamic for (a) $\eta = 0.05$ and (b) its view from the opposite side.

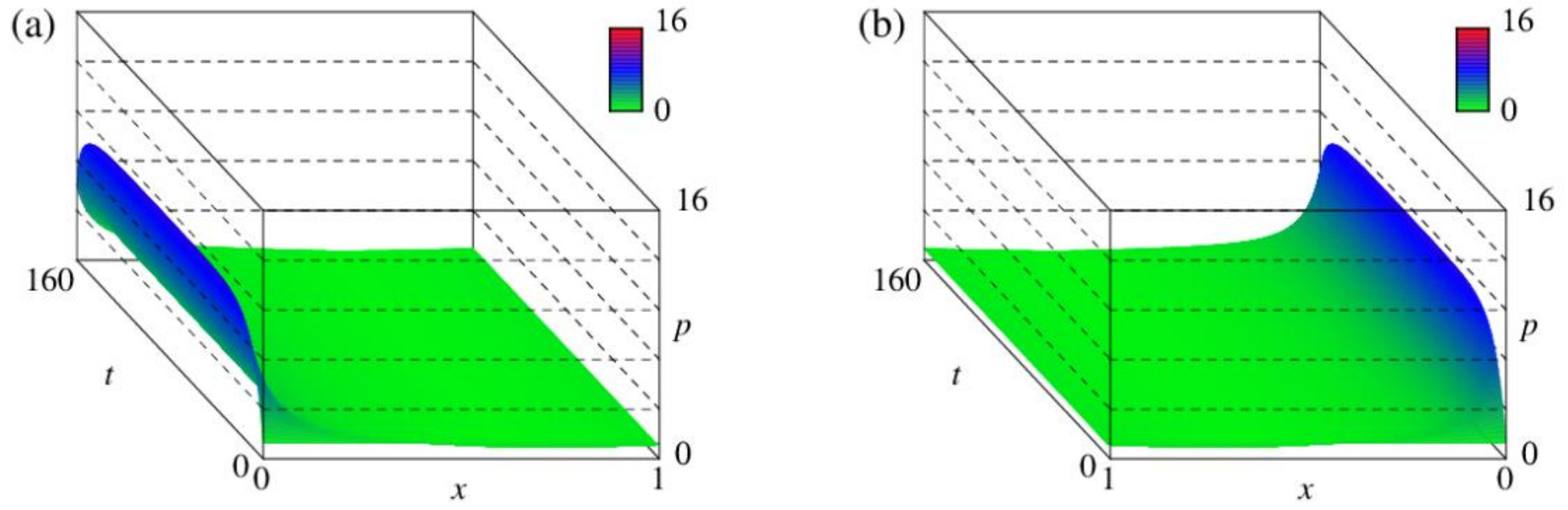

**Figure 10.** Computed PDF of GPL dynamic for (a) $\eta = 0.1$ and (b) its view from the opposite side.

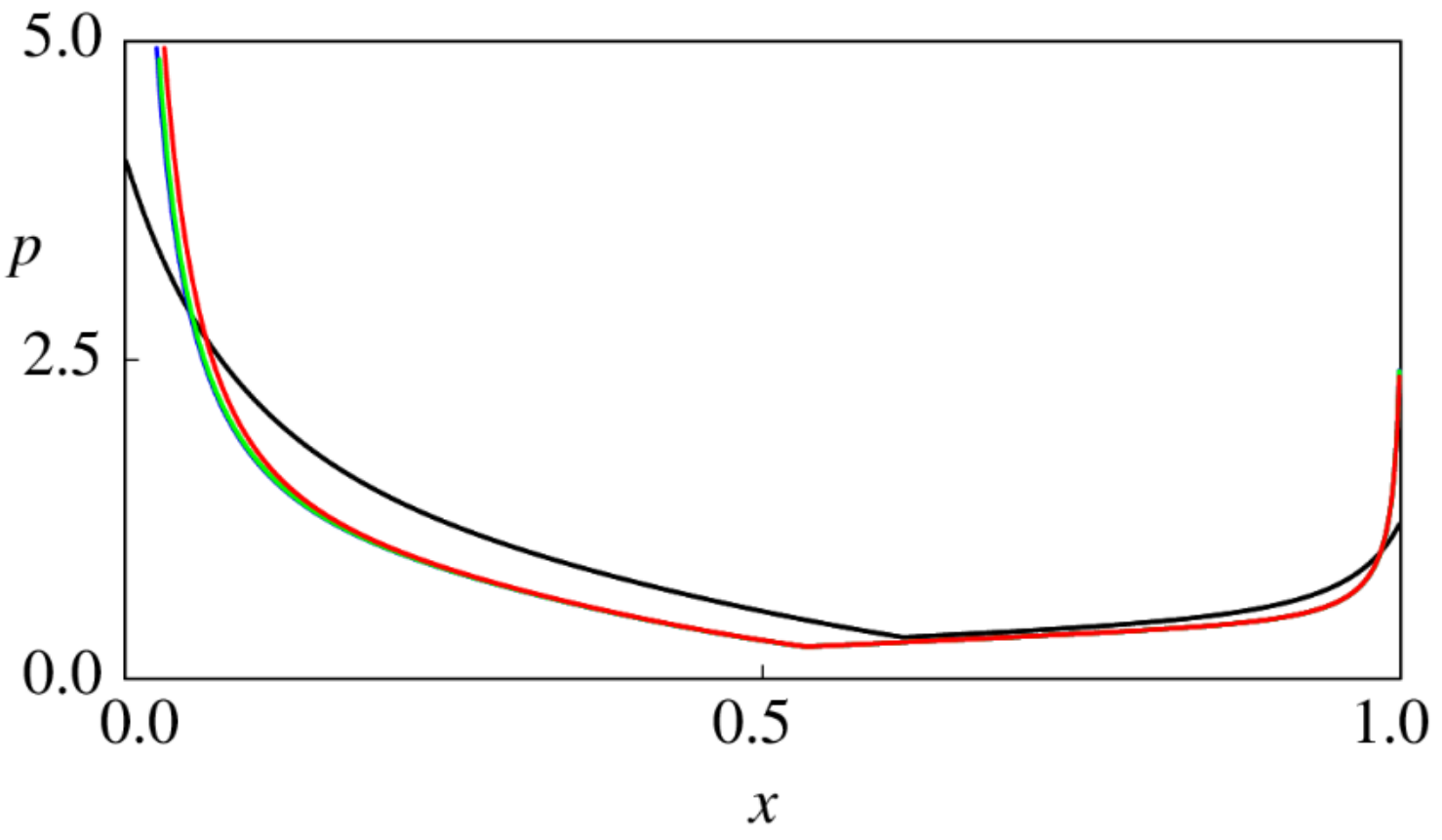

**Figure 11.** Computed (quasi-)stationary PDFs of GPL dynamic for $\eta = 0.05$ (black) and those with MFG using $T = 40$ (blue), $T = 80$ (green), $T = 160$ (red).

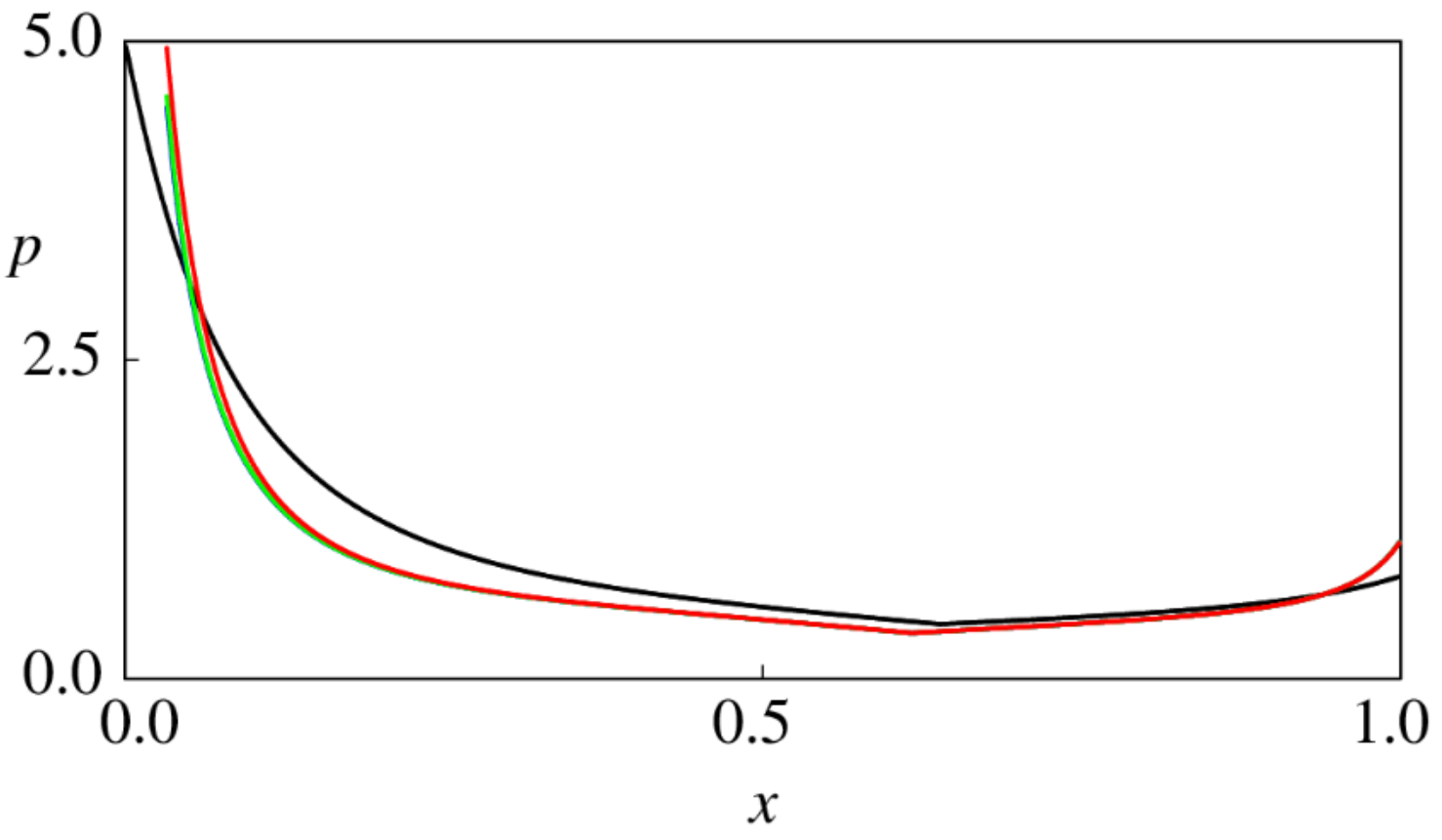

**Figure 12.** Computed (quasi-)stationary PDFs of GPL dynamic for $\eta = 0.1$ (black) and those with MFG using $T = 40$ (blue), $T = 80$ (green), $T = 160$ (red).

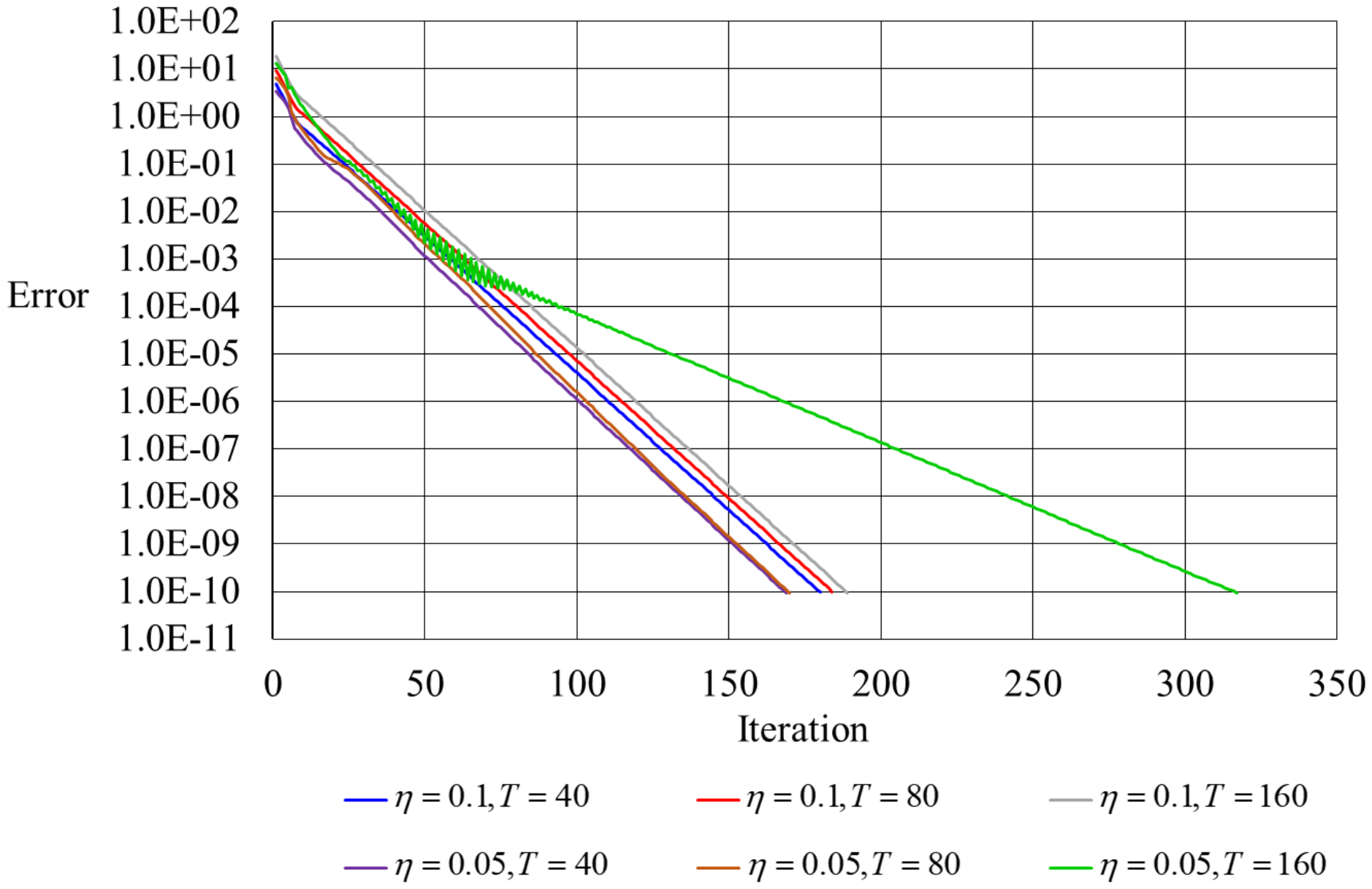

**Figure 13.** Convergence histories for MFG. Vertical axis is in logarithmic scale.

## 5. Conclusion

We proposed a GPL dynamic based on the $\kappa$-exponential function. The GPL dynamic itself is a novel dynamic, and its analysis is therefore our contribution in this study. The existence of unique solutions and approximation of the dynamic were studied. We then presented an MFG related to the GPL dynamic and discussed how the dynamic arises as a myopic version of the MFG. This version was regarded as a type of model predictive control problem. More specifically, the utility maximization encoded in the GPL dynamic was understood through a costly maximization problem where the profile of the cost determines that of the logit function. The relation of the logit dynamic to other existing approaches was also discussed. Finally, applications of the GPL dynamic and MFG to fisheries management were examined.

We focused on the pure-jump setting, although there are various applications where more generic noise processes play roles (Aurell et al., 2022; Frihi et al., 2022). The applicability and limitation of the proposed framework will then be examined. The monotonicity of the utility in our second computational example is still open, which will be addressed in the future possibly under a more generalized formulation. One can employ an approach inverse to ours, such that the logit function is designed from the objective function of an MFG if it is available from data. There will be different connections between a logit dynamic and an MFG, as recently found in the context of risk-sensitive optimization (Salhab et al., 2020). Exploring such a new connection would deepen both theoretical and engineering aspects of the evolutionary game and related research fields. Incorporating model uncertainty into the logit dynamic can be addressed based on the linkage between the evolutionary and MFGs. The other game dynamics, such as the gradient and replicator ones (Rabanal, 2017), will be examined to better understand relations between evolutionary and mean field games. Accounting for ecological dynamics (Hackney et al., 2021; Stein et al., 2023) will also be a possible next step.

## Appendices

### Appendix A. Proofs

***Proof of Lemma 1***

It suffices to show (8) when $y_1 = y_2$ because $e_\kappa(\cdot)e_\kappa(-\cdot) = 1$. We have the triangle inequality

$$\left|a(x_1, y_1, \mu_1) - a(x_2, y_1, \mu_2)\right| \le \left|a(x_1, y_1, \mu_1) - a(x_1, y_1, \mu_2)\right| + \left|a(x_1, y_1, \mu_2) - a(x_2, y_1, \mu_2)\right|. \tag{50}$$

For any $z \in \mathbb{R}$, an elementary calculation yields

$$\begin{aligned}
\frac{\mathrm{d}e_\kappa(z)}{\mathrm{d}z} &= \frac{1}{\kappa}\left(\kappa z + \sqrt{\kappa^2 z^2 + 1}\right)^{\frac{1}{\kappa}-1}\left(\kappa + \frac{\kappa^2 z}{\sqrt{\kappa^2 z^2 + 1}}\right) \\
&= e_\kappa(z)\left(\kappa z + \sqrt{\kappa^2 z^2 + 1}\right)^{-1}\left(1 + \frac{\kappa z}{\sqrt{\kappa^2 z^2 + 1}}\right) \\
&= \frac{e_\kappa(z)}{\sqrt{\kappa^2 z^2 + 1}}\left(\kappa z + \sqrt{\kappa^2 z^2 + 1}\right)^{-1}\left(\kappa z + \sqrt{\kappa^2 z^2 + 1}\right), \\
&= \frac{e_\kappa(z)}{\sqrt{\kappa^2 z^2 + 1}} \\
&\le e_\kappa(z)
\end{aligned} \tag{51}$$

and hence

$$\left|e_\kappa(z_1) - e_\kappa(z_2)\right| \le e_\kappa\left(\max\{z_1, z_2\}\right)\left|z_1 - z_2\right| \quad \text{for any } z_1, z_2 \in \mathbb{R}. \tag{52}$$

Because $e_\kappa \ge 0$, we also have

$$\begin{aligned}
&\left|a(x_1, y_1, \mu_1) - a(x_1, y_1, \mu_2)\right| \\
&= \left|\frac{1}{1 + e_\kappa\left(\frac{U(x_1, \mu_1) - U(y_1, \mu_1)}{\eta}\right)} - \frac{1}{1 + e_\kappa\left(\frac{U(x_1, \mu_2) - U(y_1, \mu_2)}{\eta}\right)}\right| \\
&= \frac{\left|e_\kappa\left(\frac{U(x_1, \mu_2) - U(y_1, \mu_2)}{\eta}\right) - e_\kappa\left(\frac{U(x_1, \mu_1) - U(y_1, \mu_1)}{\eta}\right)\right|}{\left(1 + e_\kappa\left(\frac{U(x_1, \mu_1) - U(y_1, \mu_1)}{\eta}\right)\right)\left(1 + e_\kappa\left(\frac{U(x_1, \mu_2) - U(y_1, \mu_2)}{\eta}\right)\right)} \\
&\le \left|e_\kappa\left(\frac{U(x_1, \mu_2) - U(y_1, \mu_2)}{\eta}\right) - e_\kappa\left(\frac{U(x_1, \mu_1) - U(y_1, \mu_1)}{\eta}\right)\right|
\end{aligned} \tag{53}$$

Then, combining (52) and (53) along with the boundedness $|U| \le \bar{U}$ on $\Omega \times \mathcal{M}_2$ with some $\bar{U} > 0$ yields

$$
\begin{aligned}
&\left|a(x_1,y_1,\mu_1)-a(x_1,y_1,\mu_2)\right| \\
&\le e_\kappa\left(\frac{2\bar{U}}{\eta}\right)\left|\frac{U(x_1,\mu_2)-U(y_1,\mu_2)}{\eta}-\frac{U(x_1,\mu_1)-U(y_1,\mu_1)}{\eta}\right| \\
&=\frac{1}{\eta}e_\kappa\left(\frac{2\bar{U}}{\eta}\right)\left|U(x_1,\mu_2)-U(y_1,\mu_2)-\left(U(x_1,\mu_1)-U(y_1,\mu_1)\right)\right| \\
&\le\frac{1}{\eta}e_\kappa\left(\frac{2\bar{U}}{\eta}\right)\left(\left|U(x_1,\mu_2)-U(x_1,\mu_1)\right|+\left|U(y_1,\mu_2)-U(y_1,\mu_1)\right|\right) \\
&\le\frac{2L_U}{\eta}e_\kappa\left(\frac{2\bar{U}}{\eta}\right)\left\|\mu_1-\mu_2\right\|
\end{aligned}
\quad . \tag{54}
$$

Similarly, we have

$$
\begin{aligned}
&\left|a(x_1,y_1,\mu_2)-a(x_2,y_1,\mu_2)\right| \\
&\le\left|e_\kappa\left(\frac{U(x_1,\mu_2)-U(y_1,\mu_2)}{\eta}\right)-e_\kappa\left(\frac{U(x_2,\mu_2)-U(y_1,\mu_2)}{\eta}\right)\right| \\
&\le e_\kappa\left(\frac{2\bar{U}}{\eta}\right)\left|\frac{U(x_1,\mu_2)-U(y_1,\mu_2)}{\eta}-\frac{U(x_2,\mu_2)-U(y_1,\mu_2)}{\eta}\right| \\
&=\frac{1}{\eta}e_\kappa\left(\frac{2\bar{U}}{\eta}\right)\left|U(x_1,\mu_2)-U(x_2,\mu_2)\right| \\
&\le\frac{L_U}{\eta}e_\kappa\left(\frac{2\bar{U}}{\eta}\right)\left|x_1-x_2\right|
\end{aligned}
\quad . \tag{55}
$$

From (50), (54), and (55), the desired result follows with $L_a=\frac{3L_U}{\eta}e_\kappa\left(\frac{2\bar{U}}{\eta}\right)$.

□

***Proof of Lemma 2***

Fix $M>0$, $z\in[-M,M]$, and $\kappa\in(0,1]$. We have

$$
\frac{\partial}{\partial\kappa}\ln e_\kappa(z)=\frac{1}{e_\kappa(z)}\frac{\partial e_\kappa(z)}{\partial\kappa}=-\frac{1}{\kappa}\ln e_\kappa(z)+\frac{z}{\kappa\sqrt{\kappa^2z^2+1}}, \tag{56}
$$

and hence

$$
\frac{\partial e_\kappa(z)}{\partial\kappa}=e_\kappa(z)\left(-\frac{1}{\kappa}\ln e_\kappa(z)+\frac{z}{\kappa\sqrt{\kappa^2z^2+1}}\right). \tag{57}
$$

We have

$$
\begin{aligned}
-\frac{1}{\kappa}\ln e_\kappa(z)+\frac{z}{\kappa\sqrt{\kappa^2z^2+1}}&=-\frac{1}{\kappa^2}\ln\left(\kappa z+\sqrt{\kappa^2z^2+1}\right)+\frac{\kappa z}{\kappa^2\sqrt{\kappa^2z^2+1}} \\
&=\frac{1}{\kappa^2}\left(-\ln\left(\kappa z+\sqrt{\kappa^2z^2+1}\right)+\frac{\kappa z}{\sqrt{\kappa^2z^2+1}}\right)
\end{aligned}
. \tag{58}
$$

Consequently, it follows that

$$
\frac{\partial e_\kappa(z)}{\partial\kappa}=e_\kappa(z)\frac{H(\kappa z)}{\kappa^2} \tag{59}
$$

with $H:\mathbb{R}\to\mathbb{R}$ given by

$$H(u)=-\ln\left(u+\sqrt{u^2+1}\right)+\frac{u}{\sqrt{u^2+1}} \quad \text{for any } u\in\mathbb{R}. \tag{60}$$

An elementary calculation yields

$$\frac{\mathrm{d}H(u)}{\mathrm{d}u}=-\frac{u^2}{\left(u^2+1\right)^{3/2}}\leq 0 \quad \text{for any } u\in\mathbb{R}. \tag{61}$$

This shows that $H(u)=-\frac{u^3}{3}+O\left(u^4\right)$ if $|u|$ is small. Therefore, for fixed $M,z$, it follows that

$$\lim_{\kappa\to+0}\frac{\partial e_\kappa(z)}{\partial\kappa}=\lim_{\kappa\to+0}\frac{e_\kappa(z)}{\kappa^2}H(\kappa z)=-\lim_{\kappa\to+0}\frac{e_\kappa(z)}{\kappa^2}\frac{\kappa^3 z^3}{3}=0 \tag{62}$$

Furthermore, we have $\lim_{\kappa\to+0}e_\kappa(z)=\exp(z)<+\infty$. Therefore, for each bounded $z\in\mathbb{R}$, the partial derivative $\frac{\partial e_\kappa(z)}{\partial\kappa}$ exists and is continuous for $\kappa\in[0,1]$. This proves the first statement.

Concerning the second statement, (61) shows that $H(u)\leq 0$ for $u\geq 0$; i.e., if $z>0$, then $e_\kappa(z)\leq e_0(z)=\exp(z)$. Hence, we have

$$\left|\frac{\partial e_\kappa(z)}{\partial\kappa}\right|=\frac{e_\kappa(z)}{\kappa^2}\left|H(\kappa z)\right|\leq\exp(M)\frac{\left|H(\kappa z)\right|}{\kappa^2}. \tag{63}$$

The term $\frac{\left|H(\kappa z)\right|}{\kappa^2}$ is bounded for small $\kappa\geq 0$ by (62) and the continuity of $H$. This term is also bounded for any $\kappa\in(0,1]$. Moreover, by (60) we have

$$\left|H(\kappa z)\right|\leq\ln\left(M+\sqrt{M^2+1}\right)+1 \quad \text{for any } \kappa\in[0,1] \text{ and } z\in[-M,M]. \tag{64}$$

This combined with (61)-(62) implies that there exists a constant $C'_M>0$ depending only on $M$ such that $\frac{\left|H(\kappa z)\right|}{\kappa^2}\leq C'_M$ for any $\kappa\in[0,1]$ and $|z|\leq M$. We prove the second statement by choosing $C_M=C'_M\exp(M)$.

□

***Proof of Proposition 1***

The proof for part (a) employs a two-step approach to apply existing findings on nonlinear evolution equations in Banach spaces (e.g., Proof of Theorem 3.4 of Lahkar and Riedel, 2015; Proof of Theorem 4.3 of Mendoza-Palacios and Hernández-Lerma, 2015). The first step proves a kind of Lipschitz continuity of the right-hand side of (10) with respect to $\mu\in\mathcal{M}$ ($\mu\in\mathcal{P}$ not at this stage) to show that this dynamic admits a unique solution such that $\mu(t,\cdot)\in\mathcal{M}$ for $t\in[0,+\infty)$. The second step then shows that the same solution satisfies a modified version of the dynamic, and that this solution belongs to $\mu\in\mathcal{P}$.

***First step***

We set $I:\mathcal{M}\times\mathcal{B}\to\mathbb{R}$ as

$$\begin{aligned}I(\mu,A)&=\int_A\left(\int_\Omega a(z,x,\mu)\mu(\mathrm{dz})\right)\mathrm{d}x-\int_A\left(\int_\Omega a(x,y,\mu)\mathrm{d}y\right)\mu(\mathrm{d}x)\\&:=I_1(\mu,A)-I_2(\mu,A)\end{aligned}\quad\text{for any }\mu\in\mathcal{M}\text{ and }A\in\mathcal{B}. \tag{65}$$

We show that there is a constant $L_I>0$ such that

$$\left\|I(\mu,\cdot)-I(\nu,\cdot)\right\|\le L_I\left\|\mu-\nu\right\|\quad\text{for any }\mu\in\mathcal{M}. \tag{66}$$

Fix $\mu,\nu\in\mathcal{M}_2$. We have ($\mathrm{d}x$ formally represents a Lebesgue measure)

$$\begin{aligned}\left\|I_1(\mu,\cdot)-I_1(\nu,\cdot)\right\|&=\left\|\left(\int_\Omega a(z,x,\mu)\mu(\mathrm{dz})\right)\mathrm{d}x-\left(\int_\Omega a(z,x,\nu)\nu(\mathrm{dz})\right)\mathrm{d}x\right\|\\&\le\left\|\left(\int_\Omega a(z,x,\mu)\mu(\mathrm{dz})\right)\mathrm{d}x-\left(\int_\Omega a(z,x,\mu)\nu(\mathrm{dz})\right)\mathrm{d}x\right\|\\&\quad+\left\|\left(\int_\Omega a(z,x,\mu)\nu(\mathrm{dz})\right)\mathrm{d}x-\left(\int_\Omega a(z,x,\nu)\nu(\mathrm{dz})\right)\mathrm{d}x\right\|\\&=\left\|\left(\int_\Omega a(z,x,\mu)\left\{\mu(\mathrm{dz})-\nu(\mathrm{dz})\right\}\right)\mathrm{d}x\right\|\\&\quad+\left\|\left(\int_\Omega\left\{a(z,x,\mu)-a(z,x,\nu)\right\}\nu(\mathrm{dz})\right)\mathrm{d}x\right\|\\&\le\left\|\mu-\nu\right\|+\left\|\left(\int_\Omega\left\{a(z,x,\mu)-a(z,x,\nu)\right\}\nu(\mathrm{dz})\right)\mathrm{d}x\right\|\\&\le\left\|\mu-\nu\right\|+\left\|\nu\right\|\sup_{x,z\in\Omega}\left|a(z,x,\mu)-a(z,x,\nu)\right|\\&\le\left\|\mu-\nu\right\|+2L_a\left\|\mu-\nu\right\|\end{aligned}, \tag{67}$$

where we used $|a|\le1$, $\|\nu\|\le2$, and $\int_\Omega\mathrm{d}x=1$. Hence, we obtain

$$\left\|I_1(\mu,\cdot)-I_1(\nu,\cdot)\right\|\le(2L_a+1)\left\|\mu-\nu\right\|. \tag{68}$$

Similarly, by **Lemma 1**, we have

$$\begin{aligned}\left\|I_2(\mu,\cdot)-I_2(\nu,\cdot)\right\|&=\left\|\left(\int_\Omega a(x,z,\mu)\mathrm{dz}\right)\mu(\mathrm{d}x)-\left(\int_\Omega a(x,z,\nu)\mathrm{dz}\right)\nu(\mathrm{d}x)\right\|\\&\le\left\|\left(\int_\Omega a(x,z,\mu)\mathrm{dz}\right)\mu(\mathrm{d}x)-\left(\int_\Omega a(x,z,\mu)\mathrm{dz}\right)\nu(\mathrm{d}x)\right\|\\&\quad+\left\|\left(\int_\Omega a(x,z,\mu)\mathrm{dz}\right)\nu(\mathrm{d}x)-\left(\int_\Omega a(x,z,\nu)\mathrm{dz}\right)\nu(\mathrm{d}x)\right\|,\\&\le\left\|\mu-\nu\right\|+\left\|\nu\right\|\sup_{x,z\in\Omega}\left|a(x,z,\mu)-a(x,z,\nu)\right|\\&\le(2L_a+1)\left\|\mu-\nu\right\|\end{aligned} \tag{69}$$

and hence

$$\left\|I_2(\mu,\cdot)-I_2(\nu,\cdot)\right\|\le(2L_a+1)\left\|\mu-\nu\right\|. \tag{70}$$

Consequently, with $L_I=2(2L_a+1)$ we obtain

$$\left\|I(\mu,\cdot)-I(\nu,\cdot)\right\|\le L_I\left\|\mu-\nu\right\|\quad\text{for any }\mu,\nu\in\mathcal{M}_2. \tag{71}$$

Having proven the Lipschitz continuity of $I$, we set

$$K(\mu,A)=\max\left\{2-\|\mu\|,0\right\}I(\mu,A)\quad\text{for any }\mu\in\mathcal{M}\text{ and }A\in\mathcal{B}. \tag{72}$$

Following Theorem 3.4 (p. 278) of Lahkar and Riedel (2015), by (71) and **Lemma 1**, this $K$ satisfies a Lipschitz continuity; there exists a constant $L_K > 0$ such that

$$\|K(\mu,\cdot) - K(\nu,\cdot)\| \le L_K \|\mu - \nu\| \quad \text{for any} \quad \mu,\nu \in \mathcal{M} \tag{73}$$

due to the Lipschitz continuity of $I$. More specifically, if $\|\mu\|,\|\nu\| \ge 2$ then $\|K(\mu,\cdot) - K(\nu,\cdot)\| = 0$, if $\|\mu\|,\|\nu\| \le 2$ ($\mu,\nu \in \mathcal{M}_2$) then, with the help of (68) (we used $\|I(\mu,\cdot)\| \le 4$)

$$\begin{aligned}
\|K(\mu,\cdot) - K(\nu,\cdot)\| &= \left\|\max\{2-\|\mu\|,0\} I(\mu,\cdot) - \max\{2-\|\nu\|,0\} I(\nu,\cdot)\right\| \\
&\le \left|\max\{2-\|\mu\|,0\} I(\mu,\cdot) - \max\{2-\|\nu\|,0\} I(\mu,\cdot)\right| \\
&\quad + \left|\max\{2-\|\nu\|,0\} I(\mu,\cdot) - \max\{2-\|\nu\|,0\} I(\nu,\cdot)\right| \\
&\le \|I(\mu,\cdot)\| \left\|\|\mu\| - \|\nu\|\right\| + \max\{2-\|\nu\|,0\} \|I(\mu,\cdot) - I(\nu,\cdot)\|, \\
&\le 4\left\|\|\mu\| - \|\nu\|\right\| + 2\|I(\mu,\cdot) - I(\nu,\cdot)\| \\
&\le 4\|\mu - \nu\| + 2\|I(\mu,\cdot) - I(\nu,\cdot)\| \\
&\le 2(L_I + 2)\|\mu - \nu\|
\end{aligned} \tag{74}$$

and if $\|\mu\| \le 2 \le \|\nu\|$ then (the case $\|\nu\| \le 2 \le \|\mu\|$ is similar)

$$\begin{aligned}
\|K(\mu,\cdot) - K(\nu,\cdot)\| &\le (2 - \|\mu\|)\|I(\mu,\cdot)\| \\
&\le \|I(\mu,\cdot)\|(\|\nu\| - \|\mu\|) \\
&\le \|I(\mu,\cdot)\|\|\mu - \nu\| \\
&\le 4\|\mu - \nu\|
\end{aligned}. \tag{75}$$

Then, by **Lemma 1** the Lipschitz continuity of $I$ in (74) follows with $L_K = 2(L_I + 2)$, and hence we obtain (73). Finally, the boundedness of $K$, namely $\|K(\mu,\cdot)\| \le 8$, follows analogously, by considering the cases $\|\mu\| \le 2$ and $\|\mu\| \ge 2$.

***Second step***

Based on (73), **Lemma 1**, and the generalized Picard–Lindelöf Theorem (Theorem A.3 of Lahkar et al., 2022; p. 79 of Zeidler, 1986), the following modified dynamic

$$\frac{\mathrm{d}\mu(t,A)}{\mathrm{d}t} = K(\mu,A) \quad \text{for} \quad t > 0 \tag{76}$$

subject to a prescribed initial condition $\mu(0,\cdot) = \mu_0 \in \mathcal{P}$ admits a unique solution $\mu(t,\cdot) \in \mathcal{M}$. We show that this solution is non-negative (if $\mu(0,A) \ge 0$ for any $A \in \mathcal{B}$, then $\mu(t,A) \ge 0$ ($t > 0$) for any $A \in \mathcal{B}$). To show this, fix $A \in \mathcal{B}$. If $\mu(t,A) \ge 0$ at some $t > 0$, then

$$
\begin{aligned}
\frac{\mathrm{d}\mu(t,A)}{\mathrm{d}t} &= \max\left\{2-\left\|\mu(t,\cdot)\right\|,0\right\}\left\{\int_A\left(\int_\Omega a(z,x,\mu)\mu(t,\mathrm{d}z)\right)\mathrm{d}x-\int_A\left(\int_\Omega a(x,y,\mu)\mathrm{d}y\right)\mu(t,\mathrm{d}x)\right\} \\
&\geq -\max\left\{2-\left\|\mu(t,\cdot)\right\|,0\right\}\int_A\left(\int_\Omega a(x,y,\mu)\mathrm{d}y\right)\mu(t,\mathrm{d}x) \\
&\geq -\max\left\{2-\left\|\mu(t,\cdot)\right\|,0\right\}\int_A \mu(t,\mathrm{d}x) \\
&= -2\int_A \mu(t,\mathrm{d}x) \\
&= -2\mu(t,A)
\end{aligned}
, \quad (77)
$$

and hence

$$
\frac{\mathrm{d}}{\mathrm{d}t}\left(\mu(t,A)\exp(2t)\right)\geq 0, \quad (78)
$$

which implies the non-negativity due to the arbitrariness of $A\in\mathcal{B}$ because of $\mu(0,A)\geq 0$. We then follow Proof of Theorem 4.4, Appendix A in Yegorov et al. (2020), showing the nonnegativity. We also show that the mass conservation property, namely, $\mu(t,\Omega)=1$ ($t\geq 0$) if $\mu(0,\Omega)=1$. To show this, we firstly notice that, following Proof of Theorem 4.5 of Mendoza-Palacios and Hernández-Lerma (2015) we have

$$
\frac{\mathrm{d}}{\mathrm{d}t}\left\|\mu(t,\cdot)\right\| = \frac{\mathrm{d}}{\mathrm{d}t}\int_\Omega \mu(t,\mathrm{d}x)\mathrm{d}x = \frac{\mathrm{d}\mu(t,\Omega)}{\mathrm{d}t}, \quad (79)
$$

and hence

$$
\begin{aligned}
\frac{\mathrm{d}}{\mathrm{d}t}\left\|\mu(t,\cdot)\right\| &= \frac{\mathrm{d}}{\mathrm{d}t}\int_\Omega \mu(t,\mathrm{d}x)\mathrm{d}x \\
&= \int_\Omega \max\left\{2-\left\|\mu(t,\cdot)\right\|,0\right\}\left\{\left(\int_\Omega a(z,x,\mu)\mu(t,\mathrm{d}z)\right)\mathrm{d}x-\left(\int_\Omega a(x,y,\mu)\mathrm{d}y\right)\mu(t,\mathrm{d}x)\right\} \\
&= \max\left\{2-\left\|\mu(t,\cdot)\right\|,0\right\}\left\{\int_\Omega\left(\int_\Omega a(z,x,\mu)\mu(t,\mathrm{d}z)\right)\mathrm{d}x-\left(\int_\Omega a(x,y,\mu)\mathrm{d}y\right)\mu(t,\mathrm{d}x)\right\}. \\
&= \max\left\{2-\left\|\mu(t,\cdot)\right\|,0\right\}\left\{\int_\Omega\int_\Omega a(z,x,\mu)\mu(t,\mathrm{d}z)\mathrm{d}x-\int_\Omega\int_\Omega a(x,y,\mu)\mathrm{d}y\mu(t,\mathrm{d}x)\right\} \\
&= 0
\end{aligned}
\quad (80)
$$

This shows the mass conservation property. Therefore, the unique solution to the modified dynamic (76) is non-negative and $\mu(t,\Omega)=1$ ($t\geq 0$), which is therefore a probability measure. This solution also solves the GPL dynamic, and the uniqueness follows from a classical Gronwall lemma with the help of the Lipschitz continuity of $U$ and $a$. The proof of part (b) follows from the classical estimate (e.g., the last paragraph of Proof of Theorem3.4 in Lahkar and Riedel (2015)).

□

***Proof of Proposition 2***

We use **Lemmas 1** and **2** with $M=2\bar{U}/\eta$. We write the logit function $a$ for $\kappa=\kappa_i$ as $a^{(i)}$. Fix $A\in\mathcal{B}$ and $t>0$. We have

$$
\frac{\mathrm{d}\mu^{(1)}(t,A)}{\mathrm{d}t}-\frac{\mathrm{d}\mu^{(2)}(t,A)}{\mathrm{d}t}=I^{(1)}\left(\mu^{(1)},A\right)-I^{(2)}\left(\mu^{(2)},A\right), \quad (81)
$$

where

$$I^{(i)}(\mu, A) = \int_A \left( \int_\Omega a^{(i)}(z, x, \mu) \mu(\mathrm{d}z) \right) \mathrm{d}x - \int_A \left( \int_\Omega a^{(i)}(x, y, \mu) \mathrm{d}y \right) \mu(\mathrm{d}x), \quad i = 1, 2. \tag{82}$$

Then, we obtain

$$\begin{aligned}
\frac{\mathrm{d}\mu^{(1)}(t,A)}{\mathrm{d}t} - \frac{\mathrm{d}\mu^{(2)}(t,A)}{\mathrm{d}t} &= \int_A \left( \int_\Omega a^{(1)}\left(z, x, \mu^{(1)}\right) \mu^{(1)}(t, \mathrm{d}z) \right) \mathrm{d}x - \int_A \left( \int_\Omega a^{(1)}\left(x, y, \mu^{(1)}\right) \mathrm{d}y \right) \mu^{(1)}(t, \mathrm{d}x) \\
&\quad - \int_A \left( \int_\Omega a^{(2)}\left(z, x, \mu^{(2)}\right) \mu^{(2)}(t, \mathrm{d}z) \right) \mathrm{d}x + \int_A \left( \int_\Omega a^{(2)}\left(x, y, \mu^{(2)}\right) \mathrm{d}y \right) \mu^{(2)}(t, \mathrm{d}x) \\
&= \int_A \left( \int_\Omega a^{(1)}\left(z, x, \mu^{(1)}\right) \mu^{(1)}(t, \mathrm{d}z) \right) \mathrm{d}x - \int_A \left( \int_\Omega a^{(2)}\left(z, x, \mu^{(2)}\right) \mu^{(2)}(t, \mathrm{d}z) \right) \mathrm{d}x \\
&\quad - \int_A \left( \int_\Omega a^{(1)}\left(x, y, \mu^{(1)}\right) \mathrm{d}y \right) \mu^{(1)}(t, \mathrm{d}x) + \int_A \left( \int_\Omega a^{(2)}\left(x, y, \mu^{(2)}\right) \mathrm{d}y \right) \mu^{(2)}(t, \mathrm{d}x) \\
&= \underbrace{\int_A \left( \int_\Omega a^{(1)}\left(z, x, \mu^{(1)}\right) \mu^{(1)}(t, \mathrm{d}z) \right) \mathrm{d}x - \int_A \left( \int_\Omega a^{(2)}\left(z, x, \mu^{(2)}\right) \mu^{(1)}(t, \mathrm{d}z) \right) \mathrm{d}x}_{V_1} \\
&\quad + \underbrace{\int_A \left( \int_\Omega a^{(2)}\left(z, x, \mu^{(2)}\right) \mu^{(1)}(t, \mathrm{d}z) \right) \mathrm{d}x - \int_A \left( \int_\Omega a^{(2)}\left(z, x, \mu^{(2)}\right) \mu^{(2)}(t, \mathrm{d}z) \right) \mathrm{d}x}_{V_2} \\
&\quad \underbrace{- \left\{ \int_A \left( \int_\Omega a^{(1)}\left(x, y, \mu^{(1)}\right) \mathrm{d}y \right) \mu^{(1)}(t, \mathrm{d}x) - \int_A \left( \int_\Omega a^{(2)}\left(x, y, \mu^{(1)}\right) \mathrm{d}y \right) \mu^{(2)}(t, \mathrm{d}x) \right\}}_{V_3} \\
&\quad \underbrace{- \left\{ \int_A \left( \int_\Omega a^{(2)}\left(x, y, \mu^{(1)}\right) \mathrm{d}y \right) \mu^{(2)}(t, \mathrm{d}x) - \int_A \left( \int_\Omega a^{(2)}\left(x, y, \mu^{(2)}\right) \mathrm{d}y \right) \mu^{(2)}(t, \mathrm{d}x) \right\}}_{V_4}
\end{aligned} \tag{83}$$

In the sequel, we estimate each $V_k$ ( $k = 1, 2, 3, 4$ ).

***Estimation of*** $V_1$

We have

$$\begin{aligned}
V_1 &= \int_A \left( \int_\Omega a^{(1)}\left(z, x, \mu^{(1)}\right) \mu^{(1)}(t, \mathrm{d}z) \right) \mathrm{d}x - \int_A \left( \int_\Omega a^{(2)}\left(z, x, \mu^{(2)}\right) \mu^{(1)}(t, \mathrm{d}z) \right) \mathrm{d}x \\
&= \int_A \left( \int_\Omega \left\{ a^{(1)}\left(z, x, \mu^{(1)}\right) - a^{(2)}\left(z, x, \mu^{(2)}\right) \right\} \mu^{(1)}(t, \mathrm{d}z) \right) \mathrm{d}x \\
&\le \sup_{x, z \in \Omega} \left| a^{(1)}\left(z, x, \mu^{(1)}\right) - a^{(2)}\left(z, x, \mu^{(2)}\right) \right|
\end{aligned} \tag{84}$$

For each $x, z \in \Omega$, by **Lemma 1**, we have

$$\begin{aligned}
\left| a^{(1)}\left(z, x, \mu^{(1)}\right) - a^{(2)}\left(z, x, \mu^{(2)}\right) \right| &\le \left| a^{(1)}\left(z, x, \mu^{(1)}\right) - a^{(1)}\left(z, x, \mu^{(2)}\right) \right| + \left| a^{(1)}\left(z, x, \mu^{(2)}\right) - a^{(2)}\left(z, x, \mu^{(2)}\right) \right| \\
&\le L_a \left\| \mu^{(1)}(t, \cdot) - \mu^{(2)}(t, \cdot) \right\| + \left| a^{(1)}\left(z, x, \mu^{(2)}\right) - a^{(2)}\left(z, x, \mu^{(2)}\right) \right|
\end{aligned} \tag{85}$$

The last term of (85) is evaluated as

$$\left| a^{(1)}\left(z, x, \mu^{(2)}\right) - a^{(2)}\left(z, x, \mu^{(2)}\right) \right| \le \left| e_{\kappa_1}\left( \frac{U\left(z, \mu^{(2)}\right) - U\left(x, \mu^{(2)}\right)}{\eta} \right) - e_{\kappa_2}\left( \frac{U\left(z, \mu^{(2)}\right) - U\left(x, \mu^{(2)}\right)}{\eta} \right) \right|. \tag{86}$$

Because $\left| U\left(z, \mu^{(2)}\right) - U\left(x, \mu^{(2)}\right) \right| \le L_U |x - z| \le L_U$, by taking $M = 2\eta^{-1}\bar{U}$, we apply **Lemma 2** to (86):

$$\left| a^{(1)}\left(z, x, \mu^{(2)}\right) - a^{(2)}\left(z, x, \mu^{(2)}\right) \right| \le C_M \varpi\left( |\kappa_1 - \kappa_2| \right), \tag{87}$$

where $\varpi : [0,1] \to [0, +\infty)$ is a continuous function with $\varpi(0) = 0$. Consequently, we obtain

$$V_1 \le \sup_{x,z\in\Omega}\left|a^{(1)}\left(z,x,\mu^{(1)}\right)-a^{(2)}\left(z,x,\mu^{(2)}\right)\right| \le L_a\left\|\mu^{(1)}(t,\cdot)-\mu^{(2)}(t,\cdot)\right\|+C_M\varpi\left(|\kappa_1-\kappa_2|\right). \tag{88}$$

***Estimation of*** $V_2$

Because $0\le a^{(2)}\le 1$, we have

$$\begin{aligned} V_2 &= \int_A\left(\int_\Omega a^{(2)}\left(z,x,\mu^{(2)}\right)\mu^{(1)}(t,\mathrm{dz})\right)\mathrm{d}x-\int_A\left(\int_\Omega a^{(2)}\left(z,x,\mu^{(2)}\right)\mu^{(2)}(t,\mathrm{dz})\right)\mathrm{d}x \\ &= \int_A\left(\int_\Omega a^{(2)}\left(z,x,\mu^{(2)}\right)\left\{\mu^{(1)}(t,\mathrm{dz})-\mu^{(2)}(t,\mathrm{dz})\right\}\right)\mathrm{d}x \\ &\le \left\|\mu^{(1)}(t,\cdot)-\mu^{(2)}(t,\cdot)\right\| \end{aligned}. \tag{89}$$

***Estimation of*** $V_3$

We have

$$\begin{aligned} V_3 &= -\left\{\int_A\left(\int_\Omega a^{(1)}\left(x,y,\mu^{(1)}\right)\mathrm{d}y\right)\mu^{(1)}(t,\mathrm{d}x)-\int_A\left(\int_\Omega a^{(2)}\left(x,y,\mu^{(1)}\right)\mathrm{d}y\right)\mu^{(2)}(t,\mathrm{d}x)\right\} \\ &= \int_A\left(\int_\Omega a^{(2)}\left(x,y,\mu^{(1)}\right)\mathrm{d}y\right)\mu^{(2)}(t,\mathrm{d}x)-\int_A\left(\int_\Omega a^{(1)}\left(x,y,\mu^{(1)}\right)\mathrm{d}y\right)\mu^{(1)}(t,\mathrm{d}x) \\ &= \int_A\left(\int_\Omega a^{(2)}\left(x,y,\mu^{(1)}\right)\mathrm{d}y\right)\mu^{(2)}(t,\mathrm{d}x)-\int_A\left(\int_\Omega a^{(2)}\left(x,y,\mu^{(1)}\right)\mathrm{d}y\right)\mu^{(1)}(t,\mathrm{d}x) \\ &\quad +\int_A\left(\int_\Omega a^{(2)}\left(x,y,\mu^{(1)}\right)\mathrm{d}y\right)\mu^{(1)}(t,\mathrm{d}x)-\int_A\left(\int_\Omega a^{(1)}\left(x,y,\mu^{(1)}\right)\mathrm{d}y\right)\mu^{(1)}(t,\mathrm{d}x) \\ &= \int_A\left(\int_\Omega a^{(2)}\left(x,y,\mu^{(1)}\right)\mathrm{d}y\right)\left\{\mu^{(2)}(t,\mathrm{d}x)-\mu^{(1)}(t,\mathrm{d}x)\right\} \\ &\quad +\int_A\left(\int_\Omega\left\{a^{(2)}\left(x,y,\mu^{(1)}\right)-a^{(1)}\left(x,y,\mu^{(1)}\right)\right\}\mathrm{d}y\right)\mu^{(1)}(t,\mathrm{d}x) \end{aligned}. \tag{90}$$

As in the earlier cases, we obtain

$$V_3 \le \left\|\mu^{(1)}(t,\cdot)-\mu^{(2)}(t,\cdot)\right\|+C_M\omega\left(|\kappa_1-\kappa_2|\right). \tag{91}$$

***Estimation of*** $V_4$

We have

$$\begin{aligned} V_4 &= -\left\{\int_A\left(\int_\Omega a^{(2)}\left(x,y,\mu^{(1)}\right)\mathrm{d}y\right)\mu^{(2)}(t,\mathrm{d}x)-\int_A\left(\int_\Omega a^{(2)}\left(x,y,\mu^{(2)}\right)\mathrm{d}y\right)\mu^{(2)}(t,\mathrm{d}x)\right\} \\ &= \int_A\left(\int_\Omega a^{(2)}\left(x,y,\mu^{(2)}\right)\mathrm{d}y\right)\mu^{(2)}(t,\mathrm{d}x)-\int_A\left(\int_\Omega a^{(2)}\left(x,y,\mu^{(1)}\right)\mathrm{d}y\right)\mu^{(2)}(t,\mathrm{d}x) \\ &= \int_A\left(\int_\Omega\left\{a^{(2)}\left(x,y,\mu^{(2)}\right)-a^{(2)}\left(x,y,\mu^{(1)}\right)\right\}\mathrm{d}y\right)\mu^{(2)}(t,\mathrm{d}x) \end{aligned}. \tag{92}$$

As in the earlier cases, we obtain

$$V_4 \le L_a\left\|\mu^{(1)}(t,\cdot)-\mu^{(2)}(t,\cdot)\right\|. \tag{93}$$

Consequently, by (83), (88), (89), (91), and (93), we obtain

$$\begin{aligned} \frac{\mathrm{d}\mu^{(1)}(t,A)}{\mathrm{d}t}-\frac{\mathrm{d}\mu^{(2)}(t,A)}{\mathrm{d}t} &\le L_a\left\|\mu^{(1)}(t,\cdot)-\mu^{(2)}(t,\cdot)\right\|+C_M\varpi\left(|\kappa_1-\kappa_2|\right)+\left\|\mu^{(1)}(t,\cdot)-\mu^{(2)}(t,\cdot)\right\| \\ &\quad +\left\|\mu^{(1)}(t,\cdot)-\mu^{(2)}(t,\cdot)\right\|+C_M\varpi\left(|\kappa_1-\kappa_2|\right)+L_a\left\|\mu^{(1)}(t,\cdot)-\mu^{(2)}(t,\cdot)\right\| \\ &= 2\left((L_a+1)\left\|\mu^{(1)}(t,\cdot)-\mu^{(2)}(t,\cdot)\right\|+C_M\varpi\left(|\kappa_1-\kappa_2|\right)\right) \end{aligned}. \tag{94}$$

Then, applying the classical Gronwall lemma to (94), we obtain

$$\left\|\mu^{(1)}(t,\cdot)-\mu^{(2)}(t,\cdot)\right\|\leq \tilde{C}(t)\varpi\left(\left|\kappa_1-\kappa_2\right|\right),\ \ t>0, \tag{95}$$

where $\tilde{C}(t)$ is positive and bounded at each $t\in[0,T]$, and is independent from $\kappa_1,\kappa_2$. Here, we used

$$\left|\mu^{(1)}(t,A)-\mu^{(2)}(t,A)\right|\leq \int_0^t 2\left(\left(L_a+1\right)\left\|\mu^{(1)}(s,\cdot)-\mu^{(2)}(s,\cdot)\right\|+C_M\varpi\left(\left|\kappa_1-\kappa_2\right|\right)\right) \tag{96}$$

and $\sup_A\left|\mu^{(1)}(t,A)-\mu^{(2)}(t,A)\right|=\frac{1}{2}\left\|\mu^{(1)}(t,\cdot)-\mu^{(2)}(t,\cdot)\right\|$ because both $\mu^{(1)},\mu^{(2)}$ are probability measures. The desired result directly follows from (95).

□

***Proof of Proposition 3***

We can follow Proof of Proposition 2.1 of Lahkar et al. (2023) with an adaptation to our case, especially the nonlinear term (15) not considered in this literature. We use **Lemma 1** several times in the sequel.

Fix $A\in\mathcal{B}$ and $t>0$. We have

$$\begin{aligned}\frac{\mathrm{d}\mu(t,A)}{\mathrm{d}t}-\frac{\mathrm{d}\mu_N(t,A)}{\mathrm{d}t}&=I(\mu,A)-I_N(\mu_N,A)\\&=\left(I(\mu,A)-I_N(\mu,A)\right)+\left(I_N(\mu,A)-I_N(\mu_N,A)\right)\end{aligned}. \tag{97}$$

We evaluate the last two terms of (97). For the first term in the right-hand side of (97), we have

$$\begin{aligned}I(\mu,A)-I_N(\mu,A)&=\int_A\left(\int_\Omega a(z,x,\mu)\mu(t,\mathrm{d}z)\right)\mathrm{d}x-\int_A\left(\int_\Omega a(x,y,\mu)\mathrm{d}y\right)\mu(t,\mathrm{d}x)\\&\quad-\left\{\int_A\left(\int_\Omega a_N(z,x,\mu)\mu(t,\mathrm{d}z)\right)\mathrm{d}x-\int_A\left(\int_\Omega a_N(x,y,\mu)\mathrm{d}y\right)\mu(t,\mathrm{d}x)\right\}\\&=\int_A\left(\int_\Omega\left(a(z,x,\mu)-a_N(z,x,\mu)\right)\mu(t,\mathrm{d}z)\right)\mathrm{d}x\\&\quad-\int_A\left(\int_\Omega\left(a(x,y,\mu)-a_N(x,y,\mu)\right)\mathrm{d}y\right)\mu(t,\mathrm{d}x)\end{aligned}. \tag{98}$$

Then, we obtain

$$\begin{aligned}\int_A\left(\int_\Omega\left(a(z,x,\mu)-a_N(z,x,\mu)\right)\mu(t,\mathrm{d}z)\right)\mathrm{d}x&\leq\int_A\int_\Omega\left|a(z,x,\mu)-a_N(z,x,\mu)\right|\mu(t,\mathrm{d}z)\mathrm{d}x\\&\leq\sup_{x\in\Omega}\int_\Omega\left|a(z,x,\mu)-a_N(z,x,\mu)\right|\mu(t,\mathrm{d}z)\\&\leq\sup_{x,z\in\Omega}\left|a(z,x,\mu)-a_N(z,x,\mu)\right|\end{aligned}. \tag{99}$$

We used $\int_A\mathrm{d}x\leq 1$. For any $x,z\in\Omega$, we have

$$\begin{aligned}\left|a(z,x,\mu)-a_N(z,x,\mu)\right|&\leq L_a\left|\frac{U_N(z,\mu)-U_N(x,\mu)}{\eta}-\frac{U(z,\mu)-U(x,\mu)}{\eta}\right|\\&=\frac{L_a}{\eta}\left|U_N(z,\mu)-U_N(x,\mu)-U(z,\mu)+U(x,\mu)\right|\\&\leq\frac{L_a}{\eta}\left(\left|U_N(z,\mu)-U(z,\mu)\right|+\left|U_N(x,\mu)-U(x,\mu)\right|\right)\end{aligned}. \tag{100}$$

By (17), if $z\in\Omega_i$ with some $i\in\{1,2,3,\ldots N\}$, then we have

$$
\begin{aligned}
\left|U_N(z,\mu)-U(z,\mu)\right| &= \left|U(x_{i-1/2},\mu)-U(z,\mu)\right| \\
&= \left|g\left(\sum_{j=1}^{N} f(x_{i-1/2},x_{j-1/2})\mu(t,\Omega_j)\right)-g\left(\int_\Omega f(z,y)\mu(t,\mathrm{d}y)\right)\right| \\
&\le L_g\left|\sum_{j=1}^{N} f(x_{i-1/2},x_{j-1/2})\mu(t,\Omega_j)-\int_\Omega f(z,y)\mu(t,\mathrm{d}y)\right| \\
&= L_g\left|\sum_{j=1}^{N} f(x_{i-1/2},x_{j-1/2})\mu(t,\Omega_j)-\sum_{j=1}^{N}\int_{\Omega_j} f(z,y)\mu(t,\mathrm{d}y)\right| \\
&\le L_g\sum_{j=1}^{N}\left| f(x_{i-1/2},x_{j-1/2})\mu(t,\Omega_j)-\int_{\Omega_j} f(z,y)\mu(t,\mathrm{d}y)\right| \\
&= L_g\sum_{j=1}^{N}\left|\int_{\Omega_j} f_N(x_{i-1/2},x_{j-1/2})\mu(t,\mathrm{d}y)-\int_{\Omega_j} f(z,y)\mu(t,\mathrm{d}y)\right| \\
&\le L_g\sum_{j=1}^{N}\int_{\Omega_j}\left| f_N(x_{i-1/2},x_{j-1/2})-f(z,y)\right|\mu(t,\mathrm{d}y) \\
&\le L_g\sum_{j=1}^{N}\mu(t,\Omega_j)\max_{y\in\Omega_j,z\in\Omega_i}\left| f_N(x_{i-1/2},x_{j-1/2})-f(z,y)\right| \\
&\le L_g\max_{u,v\in\Omega}\left| f_N(u,v)-f(u,v)\right|\sum_{j=1}^{N}\mu(t,\Omega_j) \\
&\le L_g\max_{u,v\in\Omega}\left| f_N(u,v)-f(u,v)\right|
\end{aligned}
\quad , \qquad (101)
$$

where we used $\sum_{j=1}^{N}\mu(t,\Omega_j)=1$. The same estimate applies to the last term of (100), yielding

$$
\sup_{x,y\in\Omega}\left|a(z,x,\mu)-a_N(z,x,\mu)\right|\le\frac{2L_aL_g}{\eta}\max_{u,v\in\Omega}\left| f_N(u,v)-f(u,v)\right| \qquad (102)
$$

and hence

$$
\int_A\left(\int_\Omega\left(a(z,x,\mu)-a_N(z,x,\mu)\right)\mu(t,\mathrm{dz})\right)\mathrm{d}x\le\frac{2L_aL_g}{\eta}\max_{u,v\in\Omega}\left| f_N(u,v)-f(u,v)\right|. \qquad (103)
$$

For the last term of (98), we have

$$
-\int_A\left(\int_\Omega\left(a(x,y,\mu)-a_N(x,y,\mu)\right)\mathrm{d}y\right)\mu(t,\mathrm{d}x)\le\int_A\int_\Omega\left|a(x,y,\mu)-a_N(x,y,\mu)\right|\mathrm{d}y\mu(t,\mathrm{d}x). \qquad (104)
$$

We again use (100)–(101), and thus (104) leads to

$$
-\int_A\left(\int_\Omega\left(a(x,y,\mu)-a_N(x,y,\mu)\right)\mathrm{d}y\right)\mu(t,\mathrm{d}x)\le\frac{2L_aL_g}{\eta}\max_{u,v\in\Omega}\left| f_N(u,v)-f(u,v)\right|. \qquad (105)
$$

Consequently, by (103) and (105), we obtain

$$
I(\mu,A)-I_N(\mu,A)\le\frac{4L_aL_g}{\eta}\max_{u,v\in\Omega}\left| f_N(u,v)-f(u,v)\right|. \qquad (106)
$$

For the second term in the right-hand side of (97), we have

$$
\begin{aligned}
I_N(\mu,A)-I_N(\mu_N,A) &= \int_A\left(\int_\Omega a_N(z,x,\mu)\mu(t,\mathrm{dz})\right)\mathrm{d}x-\int_A\left(\int_\Omega a_N(x,y,\mu)\mathrm{d}y\right)\mu(t,\mathrm{d}x) \\
&\quad -\left\{\int_A\left(\int_\Omega a_N(z,x,\mu_N)\mu_N(t,\mathrm{dz})\right)\mathrm{d}x-\int_A\left(\int_\Omega a_N(x,y,\mu_N)\mathrm{d}y\right)\mu_N(t,\mathrm{d}x)\right\} \\
&= \int_A\left(\int_\Omega a_N(z,x,\mu)\left(\mu(t,\mathrm{dz})-\mu_N(t,\mathrm{dz})\right)\right)\mathrm{d}x \\
&\quad -\left\{\int_A\left(\int_\Omega a_N(x,y,\mu)\mathrm{d}y\right)\mu(t,\mathrm{d}x)-\int_A\left(\int_\Omega a_N(x,y,\mu_N)\mathrm{d}y\right)\mu_N(t,\mathrm{d}x)\right\} \\
&= \int_A\left(\int_\Omega a_N(z,x,\mu)\left(\mu(t,\mathrm{dz})-\mu_N(t,\mathrm{dz})\right)\right)\mathrm{d}x \\
&\quad -\left\{\begin{aligned}&\int_A\left(\int_\Omega a_N(x,y,\mu)\mathrm{d}y\right)\mu(t,\mathrm{d}x)-\int_A\left(\int_\Omega a_N(x,y,\mu)\mathrm{d}y\right)\mu_N(t,\mathrm{d}x) \\ &+\int_A\left(\int_\Omega a_N(x,y,\mu)\mathrm{d}y\right)\mu_N(t,\mathrm{d}x)-\int_A\left(\int_\Omega a_N(x,y,\mu_N)\mathrm{d}y\right)\mu_N(t,\mathrm{d}x)\end{aligned}\right\} \\
&= \int_A\left(\int_\Omega a_N(z,x,\mu)\left(\mu(t,\mathrm{dz})-\mu_N(t,\mathrm{dz})\right)\right)\mathrm{d}x \\
&\quad -\left\{\int_A\left(\int_\Omega a_N(x,y,\mu)\mathrm{d}y\right)\left(\mu(t,\mathrm{d}x)-\mu_N(t,\mathrm{d}x)\right)\right\} \\
&\quad -\int_A\left(\int_\Omega\left(a_N(x,y,\mu)-a_N(x,y,\mu_N)\right)\mathrm{d}y\right)\mu_N(t,\mathrm{d}x)
\end{aligned} \quad . \tag{107}
$$

We evaluate the three terms in the last line of (107). The first term is evaluated as

$$
\begin{aligned}
\int_A\left(\int_\Omega a_N(z,x,\mu)\left(\mu(t,\mathrm{dz})-\mu_N(t,\mathrm{dz})\right)\right)\mathrm{d}x &\le \sup_{z,x\in\Omega} a_N(z,x,\mu)\left\|\mu(t,\cdot)-\mu_N(t,\cdot)\right\| \\
&\le \left\|\mu(t,\cdot)-\mu_N(t,\cdot)\right\|
\end{aligned} \quad . \tag{108}
$$

The second term is evaluated in the same way:

$$
-\left\{\int_A\left(\int_\Omega a_N(x,y,\mu)\mathrm{d}y\right)\left(\mu(t,\mathrm{d}x)-\mu_N(t,\mathrm{d}x)\right)\right\} \le \left\|\mu(t,\cdot)-\mu_N(t,\cdot)\right\| . \tag{109}
$$

The third term is evaluated as

$$
\begin{aligned}
&-\int_A\left(\int_\Omega\left(a_N(x,y,\mu)-a_N(x,y,\mu_N)\right)\mathrm{d}y\right)\mu_N(t,\mathrm{d}x) \\
&\le \int_A\left(\int_\Omega\left|a_N(x,y,\mu)-a_N(x,y,\mu_N)\right|\mathrm{d}y\right)\mu_N(t,\mathrm{d}x) \\
&\le \int_\Omega\left(\int_\Omega\left|a_N(x,y,\mu)-a_N(x,y,\mu_N)\right|\mathrm{d}y\right)\mu_N(t,\mathrm{d}x)
\end{aligned} \quad . \tag{110}
$$

As in (100), we have

$$
\left|a_N(x,y,\mu)-a_N(x,y,\mu_N)\right| \le \frac{L_a}{\eta}\left(\left|U_N(x,\mu)-U_N(x,\mu_N)\right|+\left|U_N(y,\mu)-U_N(y,\mu_N)\right|\right) . \tag{111}
$$

If $x\in\Omega_i$ with some $i\in\{1,2,3,\ldots N\}$, then we have

$$
\begin{aligned}
\left|U_N(x,\mu)-U_N(x,\mu_N)\right| &= \left|U_N(x_{i-1/2},\mu)-U_N(x_{i-1/2},\mu_N)\right| \\
&= \left|g\left(\sum_{j=1}^{n} f(x_{i-1/2},x_{j-1/2})\mu(t,\Omega_j)\right)-g\left(\sum_{j=1}^{n} f(x_{i-1/2},x_{j-1/2})\mu_N(t,\Omega_j)\right)\right| \\
&\le L_g\left|\sum_{j=1}^{n} f(x_{i-1/2},x_{j-1/2})\mu(t,\Omega_j)-\sum_{j=1}^{n} f(x_{i-1/2},x_{j-1/2})\mu_N(t,\Omega_j)\right| \\
&\le L_g\sum_{j=1}^{n} f(x_{i-1/2},x_{j-1/2})\left|\mu(t,\Omega_j)-\mu_N(t,\Omega_j)\right| \\
&\le L_g\max_{u,v\in\Omega}\left|f(u,v)\right|\sum_{j=1}^{n}\left|\mu(t,\Omega_j)-\mu_N(t,\Omega_j)\right| \\
&\le L_g\max_{u,v\in\Omega}\left|f(u,v)\right|\left\|\mu(t,\cdot)-\mu_N(t,\cdot)\right\|
\end{aligned}
. \quad (112)
$$

The same bound applies to the last term of (111), and hence,

$$
\left|a_N(x,y,\mu)-a_N(x,y,\mu_N)\right| \le \frac{2L_aL_g}{\eta}\max_{u,v\in\Omega}\left|f(u,v)\right|\left\|\mu(t,\cdot)-\mu_N(t,\cdot)\right\|. \quad (113)
$$

Substituting (113) into (110) yields

$$
\begin{aligned}
&-\int_A\left(\int_\Omega\left(a_N(x,y,\mu)-a_N(x,y,\mu_N)\right)\mathrm{d}y\right)\mu_N(t,\mathrm{d}x) \\
&\le \int_\Omega\left(\int_\Omega\frac{2L_aL_g}{\eta}\left|f(u,v)\right|\left\|\mu(t,\cdot)-\mu_N(t,\cdot)\right\|\mathrm{d}y\right)\mu_N(t,\mathrm{d}x) . \\
&= \frac{2L_aL_g}{\eta}\max_{u,v\in\Omega}\left|f(u,v)\right|\left\|\mu(t,\cdot)-\mu_N(t,\cdot)\right\|
\end{aligned}
\quad (114)
$$

Therefore, by (107), (108), (109), and (114), we obtain

$$
I_N(\mu,A)-I_N(\mu_N,A) \le 2\left\|\mu(t,\cdot)-\mu_N(t,\cdot)\right\| + \frac{2L_aL_g}{\eta}\max_{u,v\in\Omega}\left|f(u,v)\right|\left\|\mu(t,\cdot)-\mu_N(t,\cdot)\right\|. \quad (115)
$$

By (97), (98), (106), and (115), we have

$$
\begin{aligned}
&\frac{\mathrm{d}\mu(t,A)}{\mathrm{d}t}-\frac{\mathrm{d}\mu_N(t,A)}{\mathrm{d}t} \\
&= I_N(\mu,A)-I_N(\mu_N,A) \\
&\le \frac{4L_aL_g}{\eta}\max_{u,v\in\Omega}\left|f_N(u,v)-f(u,v)\right| + 2\left\|\mu(t,\cdot)-\mu_N(t,\cdot)\right\| + \frac{2L_aL_g}{\eta}\max_{u,v\in\Omega}\left|f(u,v)\right|\left\|\mu(t,\cdot)-\mu_N(t,\cdot)\right\| \\
&\le \bar{L}\left(\max_{u,v\in\Omega}\left|f_N(u,v)-f(u,v)\right| + \left\|\mu(t,\cdot)-\mu_N(t,\cdot)\right\|\right)
\end{aligned}
\quad (116)
$$

with a sufficiently large constant $\bar{L}>0$ independent from $t>0$, $N\in\mathbb{N}$, $\mu,\mu_N\in\mathcal{P}$, and $A\in\mathcal{B}$. Applying a classical Gromwell lemma combined with the arbitrariness of $A\in\mathcal{B}$ yields

$$
\left\|\mu(t,\cdot)-\mu_N(t,\cdot)\right\| \le D(t)\max_{u,v\in\Omega}\left|f_N(u,v)-f(u,v)\right|,\ t>0, \quad (117)
$$

where $D(t)$ is positive and bounded at each $t>0$. Finally, $\max_{u,v\in\Omega}\left|f_N(u,v)-f(u,v)\right|\to 0$ as $n\to+\infty$ owing to the uniform convergence of $f_N$ to $f$ on the compact set of $\Omega$ (Proof of Proposition 3.1 in Lahkar et al. (2023)). We then obtain (23).□

## Appendix B. Growth data of *P. altivelis*

Biological growth of the fish *P. altivelis* in the Hii River has been investigated by the authors. For the main document, we used data for 2020 (Yoshioka, 2023). These data were collected with the help of a union member of the HRFC. This person started harvesting the fish in the Hii River from July to October. In 2022, the arrival rate of this person at the Hii River for the harvesting was almost once every couple of days. At each harvesting day, this person harvested the fish at some location in the Hii River and measured their averaged body weights at each arrival, corresponding to the circle plots shown in **Figure B1**. The harvested fish were not released to the river, and hence different plots stand for different individuals. We fitted the logistic curve to the data using a common nonlinear least-squares method. Theoretically, a more complex growth model can be applied to the data but is not our focus in this study.

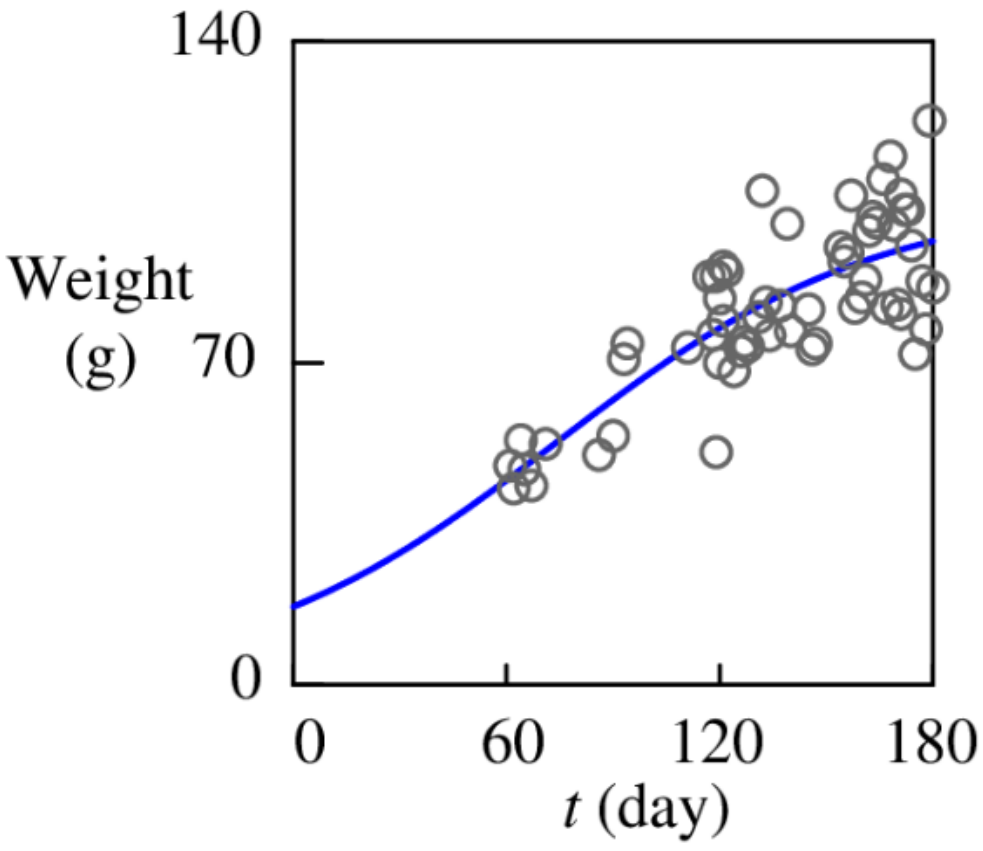

**Figure B1.** Empirical data and logistic fit of body weight of *P. altivelis* in the Hii River in 2022.

## Appendix C

### C.1 Discretization of MFG system

The MFG system to be solved has the HJB equation

$$\frac{\partial \Phi(t,x)}{\partial t} = \int_\Omega \frac{\Phi(t,x)-\Phi(t,y)}{1+e_\kappa\left(\frac{\Phi(t,x)-\Phi(t,y)}{\eta}\right)}\mathrm{d}y + \int_\Omega F\left(\frac{1}{1+e_\kappa\left(\frac{\Phi(t,x)-\Phi(t,y)}{\eta}\right)}\right)\mathrm{d}y - U\left(x,\mu(t,\cdot)\right) \tag{118}$$

for any $t < T$ and $x \in \Omega$, subject to the terminal condition $\Phi(T,\cdot) = 0$, and the FP equation

$$\frac{\mathrm{d}\mu(t,A)}{\mathrm{d}t} = \int_A \left(\int_\Omega a^*(z,x,\mu)\mu(t,\mathrm{d}z)\right)\mathrm{d}x - \int_A \left(\int_\Omega a^*(x,y,\mu)\mathrm{d}y\right)\mu(t,\mathrm{d}x) \tag{119}$$

for any $t > 0$ and $A \in \mathcal{B}$, subject to a prescribed initial condition $\mu(0,\cdot) \in \mathcal{P}$. Here, we assume $F$ of (35), and hence, $a^*(x,y,\mu) = F'^{(-1)}\left(\Phi(t,y)-\Phi(t,x)\right)$ and $G = U$. The HJB and FP equations need to be solved backward and forward in time, respectively.

The discretization method of the MFG system is based on the finite difference method in **Section 2.5**, and hence we use the notations employed there. The FP equation (119) is discretized in space as

$$\begin{aligned}&\frac{\mathrm{d}\hat{\mu}(t,\Omega_i)}{\mathrm{d}t}\\ &= \frac{1}{N}\sum_{j=1}^{N} a^*\left(x_{j-1/2},x_{i-1/2},\hat{\mu}\right)\hat{\mu}(t,\Omega_j) - \left(\frac{1}{N}\sum_{j=1}^{N} a^*\left(x_{i-1/2},x_{j-1/2},\hat{\mu}\right)\right)\mu(t,\Omega_i),\quad i=1,2,\ldots,N \text{ and } t>0,\\ &\left(=\Lambda\left(i,\hat{\mu},\hat{\Phi}\right)\right)\end{aligned} \tag{120}$$

where

$$a^*\left(x_{i-1/2},x_{j-1/2},\hat{\mu}\right) = F'^{(-1)}\left(\hat{\Phi}\left(t,x_{j-1/2}\right)-\hat{\Phi}\left(t,x_{i-1/2}\right)\right). \tag{121}$$

Here, $\hat{\Phi}\left(t,x_{i-1.2}\right)$ is the discretization of $\Phi(t,x)$ in $\Omega_i$ at time $t$. The HJB equation (118) is discretized in space as

$$\begin{aligned}&\frac{\partial \hat{\Phi}\left(t,x_{i-1/2}\right)}{\partial t}\\ &= \frac{1}{N}\sum_{j=1}^{N}\frac{\hat{\Phi}\left(t,x_{i-1/2}\right)-\hat{\Phi}\left(t,x_{j-1/2}\right)}{1+e_\kappa\left(\frac{\hat{\Phi}\left(t,x_{i-1/2}\right)-\hat{\Phi}\left(t,x_{j-1/2}\right)}{\eta}\right)}\\ &\quad+\frac{1}{N}\sum_{j=1}^{N}F\left(\frac{1}{1+e_\kappa\left(\frac{\hat{\Phi}\left(t,x_{i-1/2}\right)-\hat{\Phi}\left(t,x_{j-1/2}\right)}{\eta}\right)}\right) - U\left(x_{i-1/2},\hat{\mu}(t,\cdot)\right),\quad i=1,2,\ldots,N \text{ and } t<T.\\ &\left(=\Xi\left(i,\hat{\mu},\hat{\Phi}\right)\right)\end{aligned} \tag{122}$$

The semi-discretized system containing (120) and (122) are discretized in time using the classical forward Euler discretization. If we introduce the discrete time step $t_k = \frac{k}{M}T$ with $M \in \mathbb{N}$ ( $k = 0,1,2,...,M$ ) and $\Delta t = \frac{1}{M}T$ , by using the superscript "$(k)$" to represent the value at time $t_k$ , we set the full discretization:

$$\hat{\mu}^{(k)}(\Omega_i) = \hat{\mu}^{(k-1)}(\Omega_i) + \Lambda\left(i, \hat{\mu}^{(k-1)}, \hat{\Phi}^{(k-1)}\right)\Delta t, \quad i = 1,2,...,N \text{ and } k = 1,2,3,...,M, \tag{123}$$

and

$$\hat{\Phi}^{(k-1)}(x_{i-1/2}) = \hat{\Phi}^{(k)}(x_{i-1/2}) - \Xi\left(i, \hat{\mu}^{(k)}, \hat{\Phi}^{(k)}\right)\Delta t, \quad i = 1,2,...,N \text{ and } k = 1,2,3,...,M. \tag{124}$$

We also need to discretize the initial condition of $\mu$ , but that is not discussed here because in applications we use the uniform distribution ( $\mu \equiv 1$ in $\Omega$ ). The system (123)–(124) couples forward and backward equations and is solved using **Algorithm C1** (e.g., Barker et al., 2022). Here, the superscript "$(l)$" represents the iteration step, and $\varepsilon\left(=10^{-10}\right)$ is the error tolerance. The initial guesses of **Algorithm C1** are $\hat{\mu}_{(0)}^{(\cdot)} \equiv 0$ and $\hat{\Phi}_{(0)}^{(\cdot)} \equiv 0$ . The algorithm has a user-defined damping factor $\omega \in (0,1)$ .

***Algorithm C1***

**Step 1 (Initialization):** Prepare initial guesses $\hat{\mu}_{(0)}^{(\cdot)}$ and $\hat{\Phi}_{(0)}^{(\cdot)}$ . Initialize the iteration count $l = 0$ .

**Step2 (Computation of $\hat{\Phi}_{(l+1)}^{(\cdot)}$):** Compute (124) as $\hat{\Phi}_{(l+1)}^{(k-1)}(x_{i-1/2}) = \hat{\Phi}_{(l+1)}^{(k)}(x_{i-1/2}) - \Xi\left(i, \hat{\mu}_{(l)}^{(k)}, \hat{\Phi}_{(l+1)}^{(k)}\right)\Delta t$ at each $i = 1,2,...,N$ from $k = M$ to $k = 1$ .

**Step 3 (Computation of $\hat{\mu}_{(l+1)}^{(\cdot)}$):** Compute (123) as $\hat{\mu}_{(l+1)}^{(k)}(\Omega_i) = \hat{\mu}_{(l+1)}^{(k-1)}(\Omega_i) + \Lambda\left(i, \hat{\mu}_{(l+1)}^{(k-1)}, \hat{\Phi}_{(l+1)}^{(k-1)}\right)\Delta t$ at each $i = 1,2,...,N$ from $k = 1$ to $k = M$ .

**Step 4 (Damped update):** Update numerical solutions as $\hat{\mu}_{(l+1)}^{(\cdot)} \to \omega\hat{\mu}_{(l+1)}^{(\cdot)} + (1-\omega)\hat{\mu}_{(l)}^{(\cdot)}$ and $\hat{\Phi}_{(l+1)}^{(\cdot)} \to \omega\hat{\Phi}_{(l+1)}^{(\cdot)} + (1-\omega)\hat{\Phi}_{(l)}^{(\cdot)}$ at all $i = 1,2,...,N$ and $k = 1,2,...,M$ .

**Step 5 (Convergence check):** If $\left|\hat{\Phi}_{(l+1)}^{(k)}(x_{i-1/2}) - \hat{\Phi}_{(l)}^{(k)}(x_{i-1/2})\right|, \left|\hat{\mu}_{(l+1)}^{(k)}(x_{i-1/2}) - \hat{\mu}_{(l)}^{(k)}(x_{i-1/2})\right| \le \varepsilon$ at all $i = 1,2,...,N$ and $k = 0,1,2,...,M$ , terminate the algorithm, and the candidate $\hat{\mu}_{(l+1)}^{(\cdot)}, \hat{\Phi}_{(l+1)}^{(\cdot)}$ is regarded as the numerical solution. Otherwise, set $l \to l+1$ and go to **Step 2**.

**C.2 Convergence study**

A convergence study on the GPL dynamic and MFG counterpart is conducted for the identified model with $\eta = 0.1$ for the second example in **Section 4**. For each model, the numerical solution with the resolution $N = 2^9$ and $M = 2^4 \times 10^3$ is considered as a benchmark solution because no analytical solution has been found. Subsequently, we estimate the convergence rate of numerical solutions through a series of computations with the lower resolutions $(N,M)=(N_i,M_i)=(2^{5+i},2^i\times 10^3)$ ($i=1,2,3$). The numerical solution with $(N,M)=(N_i,M_i)$ is denoted by $\Phi_i$. Assume $(N,M)=(N_i,M_i)$. Then, the convergence rate of $\Phi_{i+1}$ is estimated for each $i=2,3$ as $\log_2\left(\frac{\|\Phi_{i+1}-\Phi_4\|_\infty}{\|\Phi_i-\Phi_4\|_\infty}\right)$ with the error $\|\Phi_i-\Phi_4\|_\infty$ being the maximum difference $|\Phi_i-\Phi_4|$ among all space–time grids using the arithmetic average; i.e., for a given $i\in\{1,2,3\}$, we set

$$|\Phi_i-\Phi_4|=\max_{\substack{j=1,2,3,\ldots,N_i\\k=1,2,3,\ldots,M_i}}\left|\Phi_i^{(k-1)}\left(x_{j-1/2}\right)-\hat{\Phi}_4\left(x_{j-1/2}\right)\Big|_{t=k/M_i}\right|, \tag{125}$$

where $\hat{\Phi}_4\left(x_{j-1/2}\right)\Big|_{t=k/M_i}$ is evaluated as the arithmetic mean of the nearest $2^{4-i}$ solution values of $\Phi_4$ defined around $x_{j-1/2}$. The arithmetic mean is obtained here because the grid points are not always overlapping between different numerical resolutions. The same method for error analysis applies to $\mu$. **Tables C.1** and **C.2** show the errors and convergence rates of the computed $\mu$ and $\Phi$ for the MFG, respectively. The first-order convergence of the obtained numerical solutions can be observed from these tables. Similar analysis has been conducted for the stationary solution to the GPL dynamic as shown in **Table C3**, again implying first-order convergence.

**C.3 Value function $\Phi$**

**Figures C1** shows the computed value functions $\Phi$ for the second application of the identified models with $T=160$ for $\eta=0.05$ and $\eta=0.1$, respectively. The figures demonstrate that the computed $\Phi$ is almost linear around the time $t=T/2$ with the slope in time approximately 0.071 and 0.53 for $\eta=0.05$ and $\eta=0.1$, respectively.

**Table C1.** Error and convergence rate of computed $\mu$ for MFG.

| $i$ | Error | Convergence rate |
|---|---|---|
| 1 | 1.52.E-01 | |
| 2 | 6.25.E-02 | 1.28.E+00 |
| 3 | 1.95.E-02 | 1.68.E+00 |

**Table C2.** Error and convergence rate of computed $\Phi$ for MFG.

| $i$ | Error | Convergence rate |
|---|---|---|
| 1 | 4.70.E-01 | |
| 2 | 1.65.E-01 | 1.51.E+00 |
| 3 | 4.89.E-02 | 1.76.E+00 |

**Table C3.** Error and convergence rate of computed $\mu$ for GPL dynamic.

| $i$ | Error | Convergence rate |
|---|---|---|
| 1 | 1.45.E-02 | |
| 2 | 6.81.E-03 | 1.09.E+00 |
| 3 | 2.38.E-03 | 1.52.E+00 |

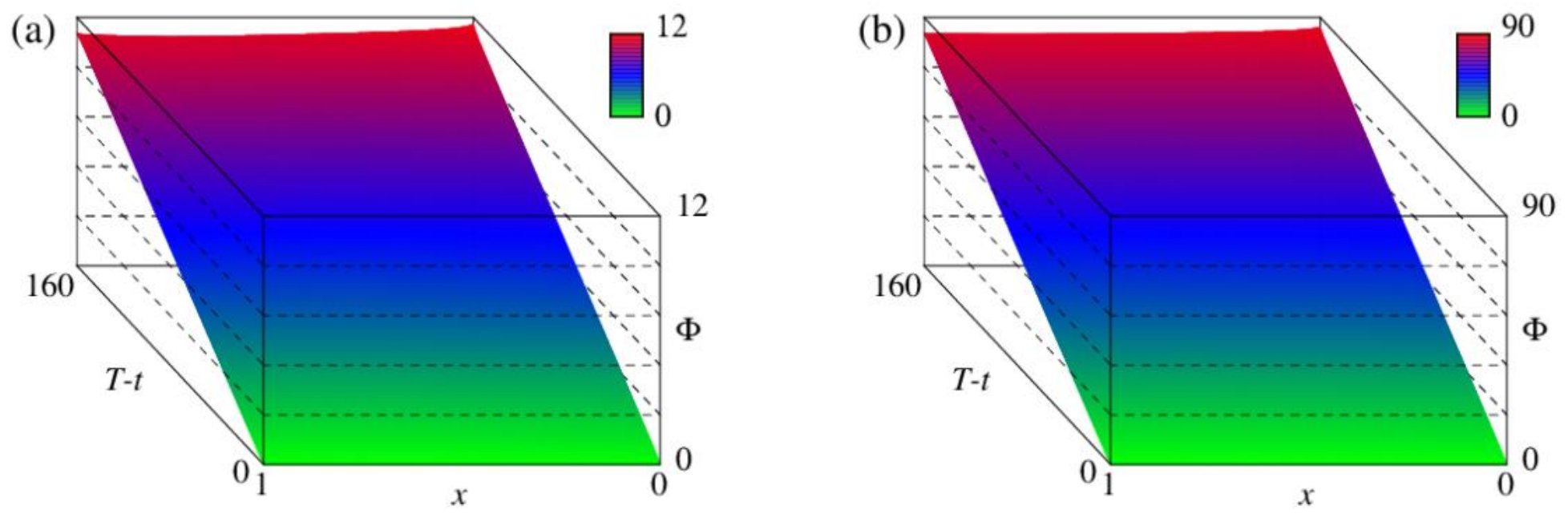

**Figure C1.** Computed value functions $\Phi$ for second application of identified models with $T = 160$: (a) $\eta = 0.05$ and (b) $\eta = 0.1$

**Declarations**

**Ethical Approval** N.A.
**Competing interests:** The authors declare no competing interests.
**Acknowledgements:** The authors thank the union members of the Hii River Fisheries Cooperative for their valuable support in our field survey.
**Funding:** Japan Society for the Promotion of Science (22K14441 and 22H02456).
**Availability of data and material:** Data are available upon reasonable request to the corresponding author.
**Declaration:** No AI technologies were used for writing the manuscript.
**Author contributions:**
*Hidekazu Yoshioka*: Conceptualization, Methodology, Software, Formal analysis, Data acquisition, Data curation, Visualization, Field survey, Writing–original draft preparation, Writing–review and editing, Supervision, Project administration, Funding acquisition.
*Motoh Tsujimura*: Formal analysis, Writing–review and editing.

**References**

[1] Achdou Y, Capuzzo-Dolcetta I (2010) Mean field games: numerical methods. SIAM Journal on Numerical Analysis, 48(3), 1136-1162. https://doi.org/10.1137/090758477
[2] Ahmadi Forushani Z, Safavi HR, Kerachian R, Golmohammadi MH (2023) A signaling game model for evaluating water allocation competitiveness with information asymmetry environment, case study: Zayandehrud River Basin, Iran. Environment, Development and Sustainability, published online. https://doi.org/10.1007/s10668-023-03989-1
[3] Aiba K (2015) Waiting times in evolutionary dynamics with time-decreasing noise. International Journal of Game Theory, 44, 499-514. https://doi.org/10.1007/s00182-014-0442-x
[4] Apicella A, Donnarumma F, Isgrò F, Prevete R (2021) A survey on modern trainable activation functions. Neural Networks, 138, 14-32. https://doi.org/10.1016/j.neunet.2021.01.026
[5] Aurell A, Carmona R, Dayanikli G, Lauriere M (2022) Optimal incentives to mitigate epidemics: a Stackelberg mean field game approach. SIAM Journal on Control and Optimization, 60(2), S294-S322. https://doi.org/10.1137/20M1377862
[6] Barker M (2019) From mean field games to the best reply strategy in a stochastic framework. Journal of Dynamics and Games, 6(4), 291-314. 10.3934/jdg.2019020
[7] Barker M, Degond P, Wolfram MT (2022) Comparing the best-reply strategy and mean-field games: The stationary case. European Journal of Applied Mathematics, 33(1), 79-110. https://doi.org/10.1017/S0956792520000376
[8] Barreiro-Gomez J, Poveda JI (2023) Recurrent Neural Network ODE Output for Classification Problems Follows the Replicator Dynamics. IEEE Control Systems Letters, 7, 3783-3788. 10.1109/LCSYS.2023.3341096
[9] Bayraktar E, Cecchin A, Chakraborty P (2023) Mean field control and finite agent approximation for regime-switching jump diffusions. Applied Mathematics & Optimization, 88(2), 36. https://doi.org/10.1007/s00245-023-10015-3
[10] Bect J (2010) A unifying formulation of the Fokker–Planck–Kolmogorov equation for general stochastic hybrid systems. Nonlinear Analysis: Hybrid Systems, 4(2), 357-370. https://doi.org/10.1016/j.nahs.2009.07.008
[11] Benjamin C, Arjun Krishnan UM, (2023) Nash equilibrium mapping vs. Hamiltonian dynamics vs. Darwinian evolution for some social dilemma games in the thermodynamic limit. The European Physical Journal B, 96(7), 105. https://doi.org/10.1140/epjb/s10051-023-00573-4
[12] Biancardi M, Iannucci G, Villani G (2023) Groundwater management and illegality in a differential-evolutionary framework. Computational Management Science, 20(1), 16. https://doi.org/10.1007/s10287-023-00449-z
[13] Carmona R, Cormier Q, Soner HM (2023) Synchronization in a Kuramoto mean field game. Communications in Partial Differential Equations, 48(9), 1214-1244. https://doi.org/10.1080/03605302.2023.2264611
[14] Cartea Á, Jaimungal S, Qin Z (2019) Speculative trading of electricity contracts in interconnected locations. Energy Economics, 79, 3-20. https://doi.org/10.1016/j.eneco.2018.11.019
[15] Cecchin A, Pra PD, Fischer M, Pelino G (2019) On the convergence problem in mean field games: a two state model without uniqueness. SIAM Journal on Control and Optimization, 57(4), 2443-2466. https://doi.org/10.1137/18M1222454
[16] Certório J, Martins NC, La RJ (2022) Epidemic population games with nonnegligible disease death rate. IEEE Control Systems Letters, 6, 3229-3234. DOI: 10.1109/LCSYS.2022.3183477
[17] Chen F, Jiang Y, Liu Z, Lin R, Yang W (2023) Framework system of marine sustainable development assessment based on systematic review. Marine Policy, 154, 105689. https://doi.org/10.1016/j.marpol.2023.105689
[18] Cheung MW (2016) Imitative dynamics for games with continuous strategy space. Games and Economic Behavior, 99, 206-223. https://doi.org/10.1016/j.geb.2016.08.003
[19] Choudhury KD, Aydinyan T (2023) Stochastic replicator dynamics: A theoretical analysis and an experimental assessment. Games and Economic Behavior, 142, 851-865. https://doi.org/10.1016/j.geb.2023.10.004
[20] Cohen A, Zell E (2023) Analysis of the finite-state ergodic master equation. Applied Mathematics & Optimization, 87(3), 40. https://doi.org/10.1007/s00245-022-09954-0

[21] Dai Pra P, Sartori E, Tolotti M (2023) Polarization and coherence in mean field games driven by private and social utility. Journal of Optimization Theory and Applications, 198, 49-85. https://doi.org/10.1007/s10957-023-02233-0
[22] Delarue F, Tchuendom RF (2020) Selection of equilibria in a linear quadratic mean-field game. Stochastic Processes and their Applications, 130(2), 1000-1040. https://doi.org/10.1016/j.spa.2019.04.005
[23] Dombi J, Jónás T (2022) The generalized sigmoid function and its connection with logical operators. International Journal of Approximate Reasoning, 143, 121-138. https://doi.org/10.1016/j.ijar.2022.01.006
[24] Drechsler M (2023a) Ecological and economic trade-offs between amount and spatial aggregation of conservation and the cost-effective design of coordination incentives. Ecological Economics, 213, 107948. https://doi.org/10.1016/j.ecolecon.2023.107948
[25] Drechsler M (2023b) Insights from Ising models of land-use under economic coordination incentives. Physica A: Statistical Mechanics and its Applications, 625, 128987. https://doi.org/10.1016/j.physa.2023.128987
[26] Dupret JL, Hainaut D (2023) A subdiffusive stochastic volatility jump model. Quantitative Finance, 23(6), 979-1002. https://doi.org/10.1080/14697688.2023.2199959
[27] Friedman D, Ostrov DN (2013) Evolutionary dynamics over continuous action spaces for population games that arise from symmetric two-player games. Journal of Economic Theory, 148(2), 743-777. https://doi.org/10.1016/j.jet.2012.07.004
[28] Frihi ZEO, Choutri SE, Barreiro-Gomez J, Tembine H (2022) Hierarchical mean-field type control of price dynamics for electricity in smart grid. Journal of Systems Science and Complexity, 35(1), 1-17. https://doi.org/10.1007/s11424-021-0176-3
[29] Fu R, Liu J (2023) Revenue sources of natural resources rents and its impact on sustainable development: Evidence from global data. Resources Policy, 80, 103226. https://doi.org/10.1016/j.resourpol.2022.103226
[30] Huang X (2023) Transboundary watershed pollution control analysis for pollution abatement and ecological compensation. Environmental Science and Pollution Research, 30(15), 44025-44042. https://doi.org/10.1007/s11356-023-25177-4
[31] Gao L, Pan Q, He M (2024) Impact of dynamic compensation with resource feedback on the common pool resource game. Chaos, Solitons & Fractals, 180, 114545. https://doi.org/10.1016/j.chaos.2024.114545
[32] Ghilli D, Ricci C, Zanco G (2023). A mean field game model for COVID-19 with human capital accumulation. Economic Theory, 77, 533-560. https://doi.org/10.1007/s00199-023-01505-0
[33] Gomes DA, Mohr J, Souza RR (2013) Continuous time finite state mean field games. Applied Mathematics & Optimization, 68(1), 99-143. https://doi.org/10.1007/s00245-013-9202-8
[34] Gomes DA, Saúde J (2021) Numerical methods for finite-state mean-field games satisfying a monotonicity condition. Applied Mathematics & Optimization, 83(1), 51-82. https://doi.org/10.1007/s00245-018-9510-0
[35] Hackney M, James A, Plank MJ (2021) Cooperative and non-cooperative behaviour in the exploitation of a common renewable resource with environmental stochasticity. Applied Mathematical Modelling, 89, 1041-1054. https://doi.org/10.1016/j.apm.2020.06.079
[36] Hata M, Otake T (2019) Downstream migration and mortality of larval Ayu (*Plecoglossus altivelis altivelis*) on the Sanriku coast, northern Japan. Environmental Biology of Fishes, 102, 1311-1325. https://doi.org/10.1007/s10641-019-00909-z
[37] Hii River Fishery Cooperative (HRFC) (2023) https://www.hiikawafish.jp/chet (in Japanese). Last accessed on December 19, 2023.
[38] Inui R et al. (2021) Spatiotemporal changes of the environmental DNA concentrations of amphidromous fish *Plecoglossus altivelis altivelis* in the spawning grounds in the Takatsu River, western Japan. Frontiers in Ecology and Evolution, 9, 622149. https://doi.org/10.1021/acs.est.2c01904
[39] Kaniadakis G et al. (2020) The κ-statistics approach to epidemiology. Scientific reports, 10(1), 19949. https://doi.org/10.1038/s41598-020-76673-3
[40] Keliger D, Horvath I, Takacs B (2022) Local-density dependent Markov processes on graphons with epidemiological applications. Stochastic Processes and their Applications, 148, 324-352. https://doi.org/10.1016/j.spa.2022.03.001
[41] Komaee A (2021) An inverse optimal approach to design of feedback control: Exploring analytical solutions for the Hamilton–Jacobi–Bellman equation. Optimal Control Applications and Methods, 42(2), 469-485. https://doi.org/10.1002/oca.2686
[42] Kroell E, Pesenti SM, Jaimungal S (2023) Stressing Dynamic Loss Models. Insurance: Mathematics and Economics, 114, 56-78. https://doi.org/10.1016/j.insmatheco.2023.11.002
[43] Lahkar R, Riedel F (2015) The logit dynamic for games with continuous strategy sets. Games and Economic Behavior, 91, 268-282. https://doi.org/10.1016/j.geb.2015.03.009
[44] Lahkar R, Mukherjee S, Roy S (2022) Generalized perturbed best response dynamics with a continuum of strategies. Journal of Economic Theory, 200, 105398. https://doi.org/10.1016/j.jet.2021.105398
[45] Lahkar R, Mukherjee S, Roy S (2023) The logit dynamic in supermodular games with a continuum of strategies: A deterministic approximation approach. Games and Economic Behavior, 139, 133-160. https://doi.org/10.1016/j.geb.2023.02.003
[46] Lai C, Li R, Gao X (2024) Bank competition with technological innovation based on evolutionary games. International Review of Economics & Finance, 89, 742-759. https://doi.org/10.1016/j.iref.2023.07.055

[47] Lasry JM, Lions PL (2007) Mean field games. Japanese journal of mathematics, 2(1), 229-260. https://doi.org/10.1007/s11537-007-0657-8

[48] Laurière M, Tangpi L (2022) Convergence of large population games to mean field games with interaction through the controls. SIAM Journal on Mathematical Analysis, 54(3), 3535-3574. https://doi.org/10.1137/22M1469328

[49] Leon V, Etesami SR, Nagi R (2023) Limited-Trust in Diffusion of Competing Alternatives Over Social Networks. IEEE Transactions on Network Science and Engineering, 11(1), 1320-1336. 10.1109/TNSE.2023.3322132

[50] Leonardos S, Piliouras G (2022) Exploration-exploitation in multi-agent learning: Catastrophe theory meets game theory. Artificial Intelligence, 304, 103653. https://doi.org/10.1016/j.artint.2021.103653

[51] Lin R et al. (2022) Optimal scheduling management of the parking lot and decentralized charging of electric vehicles based on Mean Field Game. Applied Energy, 328, 120198. https://doi.org/10.1016/j.apenergy.2022.120198

[52] Magnus R (2023) Essential Ordinary Differential Equations. Springer, Cham.

[53] Mendoza-Palacios S, Hernández-Lerma O (2015) Evolutionary dynamics on measurable strategy spaces: Asymmetric games. Journal of Differential Equations, 259(11), 5709-5733. https://doi.org/10.1016/j.jde.2015.07.005

[54] Mendoza-Palacios S, Hernández-Lerma O (2020) The replicator dynamics for games in metric spaces: finite approximations. In Advances in Dynamic Games: Games of Conflict, Evolutionary Games, Economic Games, and Games Involving Common Interest (pp. 163-186). Birkhäuser, Cham. https://doi.org/10.1007/978-3-030-56534-3_7

[55] Ministry of Agriculture, Forestry and Fisheries (2023) https://www.maff.go.jp/j/tokei/kouhyou/naisui_gyosei/. Last accessed on December 6, 2023.

[56] Nagayama S, Fujii R, Harada M, Sueyoshi M (2023) Low water temperature and increased discharge trigger downstream spawning migration of ayu *Plecoglossus altivelis*. Fisheries Science, 89, 463-475. https://doi.org/10.1007/s12562-023-01694-6

[57] Nakayama S (2013) $q$-Generalized Logit Route Choice an Network Equilibrium Model. Procedia-Social and Behavioral Sciences, 80, 753-763. https://doi.org/10.1016/j.sbspro.2013.05.040

[58] Neumann BA (2023a) Nonlinear Markov chains with finite state space: Invariant distributions and long-term behaviour. Journal of Applied Probability, 60(1), 30-44. https://doi.org/10.1017/jpr.2022.23

[59] Neumann BA (2023b) A myopic adjustment process for mean field games with finite state and action space. International Journal of Game Theory, 53, 159-195. https://doi.org/10.1007/s00182-023-00866-z

[60] Øksendal B, Sulem A (2019) Applied Stochastic Control of Jump Diffusions. Springer, Cham.

[61] Parise F, Ozdaglar A (2019) Graphon games. In Proceedings of the 2019 ACM Conference on Economics and Computation (pp. 457-458). https://doi.org/10.1145/3328526.3329638

[62] Perkins S. Leslie DS (2014) Stochastic fictitious play with continuous action sets. Journal of Economic Theory, 152, 179-213. https://doi.org/10.1016/j.jet.2014.04.008

[63] Petrov I (2023) Structural Interventions in Linear Best-Response Games on Random Graphs. IFAC-PapersOnLine, 56(2), 2830-2833. https://doi.org/10.1016/j.ifacol.2023.10.1396

[64] Rabanal JP (2017) On the evolution of continuous types under replicator and gradient dynamics: two examples. Dynamic Games and Applications, 7, 76-92. https://doi.org/10.1007/s13235-015-0164-0

[65] Rettieva A (2023) Cooperation maintenance in dynamic discrete-time multicriteria games with application to bioresource management problem. Journal of Computational and Applied Mathematics, 115699. https://doi.org/10.1016/j.cam.2023.115699

[66] Roy A, Singh C, Narahari Y (2023) Recent advances in modeling and control of epidemics using a mean field approach. Sādhanā, 48(4), 1-20. https://doi.org/10.1007/s12046-023-02268-z

[67] Salhab R, Malhamé RP, Le Ny J (2020) Collective stochastic discrete choice problems: A min-LQG dynamic game formulation. IEEE Transactions on Automatic Control, 65(8), 3302-3316. 10.1109/TAC.2019.2941443

[68] Sandholm WH (2020) Evolutionary game theory. Complex Social and Behavioral Systems: Game Theory and Agent-Based Models, 573-608. https://doi.org/10.1007/978-1-0716-0368-0_188

[69] Sandholm WH, Staudigl M (2018) Sample path large deviations for stochastic evolutionary game dynamics. Mathematics of Operations Research, 43(4), 1348-1377. https://doi.org/10.1287/moor.2017.0908

[70] Sartzetakis E, Xepapadeas A, Yannacopoulos AN (2023) Environmental regulation with preferences for social status. Ecological Economics, 209, 107834. https://doi.org/10.1016/j.ecolecon.2023.107834

[71] Stein A et al. (2023) Stackelberg evolutionary game theory: how to manage evolving systems. Philosophical Transactions of the Royal Society B, 378(1876), 20210495. https://doi.org/10.1098/rstb.2021.0495

[72] Takashina N, Cheung H, Miyazawa M (2023) Spread the word: Sharing information on social media can stabilize conservation funding and improve ecological outcomes. Conservation Science and Practice, 5(5), e12857. https://doi.org/10.1111/csp2.12857

[73] Tan S, Fang Z, Wang Y, Lü J (2023) Consensus-based multipopulation game dynamics for distributed nash equilibria seeking and optimization. IEEE Transactions on Systems, Man, and Cybernetics: Systems, 53(2), 813-823. 10.1109/TSMC.2022.3188266

[74] Thon FM, Müller C, Wittmann MJ (2024) The evolution of chemodiversity in plants-From verbal to quantitative models. Ecology Letters, 27(2), e14365. https://doi.org/10.1111/ele.14365

[75] Ullmo D, Swiecicki I, Gobron T (2019) Quadratic mean field games. Physics Reports, 799, 1-35. https://doi.org/10.1016/j.physrep.2019.01.001

[76] Ulucak R, Baloch MA (2023) An empirical approach to the nexus between natural resources and environmental pollution: Do economic policy and environmental-related technologies make any difference?. Resources Policy, 81, 103361. https://doi.org/10.1016/j.resourpol.2023.103361

[77] Wang X, Shen T, Fan E, Shen C (2023) Game theory-based stakeholder analysis of marine nature reserves and its case studies in Guangdong Province, China. Journal for Nature Conservation, 71, 126322. https://doi.org/10.1016/j.jnc.2022.126322

[78] Wang Y, Xin B (2023) Stability of a mean field game of production adjustment. Asian Journal of Control, 26(2), 717-727. https://doi.org/10.1002/asjc.3199

[79] Yano R, Kuroda H (2023) Mean-field game analysis of crowd evacuation using the Cristiani-Santo-Menci method. Physical Review E, 108(1), 014119. https://doi.org/10.1103/PhysRevE.108.014119

[80] Yegorov I, Novozhilov AS, Bratus AS (2020) Open quasispecies models: Stability, optimization, and distributed extension. Journal of Mathematical Analysis and Applications, 481(2), 123477. https://doi.org/10.1016/j.jmaa.2019.123477

[81] Yu J, Lai R, Li W, Osher S (2023) Computational mean-field games on manifolds. Journal of Computational Physics, 484, 112070. https://doi.org/10.1016/j.jcp.2023.112070

[82] Yoshimura M, Matsuura T, Sugimura K (2021) Attitudes to forest conditions and fishing activities in the mountain area in Japan. Fisheries Research, 244, 106125. https://doi.org/10.1016/j.fishres.2021.106125

[83] Yoshioka H (2023) Optimal harvesting policy for biological resources with uncertain heterogeneity for application in fisheries management, Natural Resource Modeling, e12394. https://doi.org/10.1111/nrm.12394

[84] Yoshioka H (2024) Generalized logit dynamics based on rational logit functions, Dynamic Games and Applications. https://doi.org/10.1007/s13235-023-00551-6

[85] Yoshioka H, Hamagami K, Tomobe H (2023) A Non-local Fokker-Planck Equation with Application to Probabilistic Evaluation of Sediment Replenishment Projects. Methodology and Computing in Applied Probability, 25(1), 34. https://doi.org/10.1007/s11009-023-10006-5

[86] Yoshioka H, Tsujimura M, Yoshioka Y (2024) A rational logit dynamic for decision-making under uncertainty: well-posedness, vanishing-noise limit, and numerical approximation. Accepted as a full paper to appear in the ICCS 2024 – Málaga, Spain. Preprint version is available at https://arxiv.org/abs/2402.13453

[87] Yoshioka H, Yaegashi Y (2016) Finding the optimal opening time of harvesting farmed fishery resources. Pacific journal of mathematics for industry, 8(1), 1-6. https://doi.org/10.1186/s40736-016-0025-9

[88] Zeidler E (1986) Nonlinear Functional Analysis and its Applications I: Fixed-Point Theorems, Springer, New York.

[89] Zhang Y, Wang X, Huang Q, Duan J, Li T (2020) Numerical analysis and applications of Fokker-Planck equations for stochastic dynamical systems with multiplicative α-stable noises. Applied Mathematical Modelling, 87, 711-730. https://doi.org/10.1016/j.apm.2020.06.031

[90] Zusai D (2023) Evolutionary dynamics in heterogeneous populations: a general framework for an arbitrary type distribution. International Journal of Game Theory, 52, 1215-1260. https://doi.org/10.1007/s00182-023-00867-y